\documentclass[11pt,twoside]{article} 
\usepackage{amsmath} 
\usepackage{amsthm} 
\usepackage{amssymb} 
\usepackage{amsfonts} 
\usepackage{epsfig} 
 
\input xypic

\let\MathTimeFonts= N 
\let\FormalScriptFont= N

\if Y\MathTimeFonts% 
    \usepackage{mathtime}  
\SetMathAlphabet{\mathit}{normal}{\encodingdefault}{\rmdefault}{m}{it}% 
\fi 
 
\if Y\FormalScriptFont% 
   \usepackage{calrsfs} 
\fi 
 
\makeatletter 
 
\newtheorem{The}{Theorem}[section] 
\newtheorem{Cor}[The]{Corollary} 
\newtheorem{Pro}[The]{Proposition} 
 
\newtheorem{Rem}[The]{Remark} 
\newtheorem{Lem}[The]{Lemma} 
\newtheorem{Def}[The]{Definition}

\renewcommand{\theenumi}{\alph{enumi}}

\renewcommand{\p@enumi}{} 
\renewcommand{\p@enumii}{(\theenumi)}

\def\proof{\vspace{2ex}\noindent{\bf Proof.} } 
\def\endproof{\relax\ifmmode\expandafter\endproofmath\else 
\unskip\nobreak\hfil\penalty50\hskip.75em\hbox{}\nobreak\hfil\bull 
{\parfillskip= 0pt \finalhyphendemerits= 0 \bigbreak}\fi} 
\def\endproofmath$${\eqno\bull$$\bigbreak} 
\def\bull{\vbox{\hrule\hbox{\vrule\kern3pt\vbox{\kern6pt}\kern3pt\vrule} 
\hrule}}  
 
\def\cancel#1#2{\ooalign{$\hfil#1\mkern1mu/\hfil$\crcr$#1#2$}} 
 
\def\dirac{\mathpalette\cancel\partial} 
 
\newcommand{\ba}{\begin{eqnarray}} 
\newcommand{\na}{\end{eqnarray}} 
\newcommand{\ban}{\begin{eqnarray*}} 
\newcommand{\nan}{\end{eqnarray*}} 
\newcommand{\spinc}{\mathrm{Spin}^c} 
 
\newcommand{\scr}{\mathcal} 
\newcommand{\C}{{\cal C}}
\newcommand{\CC}{\mathbb{C}} 
 
\newcommand{\R}{\mathbb{R}} 
\newcommand{\Z}{\mathbb{Z}}

\newcommand{\G}{\scr{G}} 
\newcommand{\B}{\scr{B}} 
\newcommand{\M}{\scr{M}}

\newcommand{\la}{\langle} 
\newcommand{\ra}{\rangle}

\newcommand{\nab}{\nabla} 
 
\newcommand{\s}{\mathfrak{s}}

\newcommand{\p}{\mathfrak{p}} 
\renewcommand{\d}{\partial} 
 
\newcommand{\disp}{\displaystyle}

\title{The geometric triangle for 3-dimensional Seiberg-Witten monopoles} 
\author{A.L. Carey, M. Marcolli, B.L. Wang} 
\date{}

\begin{document} 
 
\maketitle 
 
\tableofcontents

\section{Introduction} 
 
This paper is the first part of a program aimed at 
a better understanding of  how the recently defined 
Seiberg-Witten-Floer homology for any closed 3-manifold 
$Y$ with a $\spinc$ structure $\s$ \cite{CW}, \cite{K}, \cite{Mar2}, 
\cite{MW}, \cite{RW1} behaves under surgery. 
The non-equivariant Seiberg-Witten-Floer homology is constructed 
from the chain complex generated by the irreducible critical points 
of the perturbed Chern-Simons-Dirac functional  
on the space of $L^2_1$-configurations modulo the  
action of $L^2_2$-gauge transformations, the differential 
is defined by counting the gradient flow lines connecting the critical 
points of relative index one. These critical points are the  
equivalence classes of solutions to the Seiberg-Witten equations 
on $(Y, \s)$ modulo gauge transformations. 
The gradient flowlines are the equivalence classes of solutions to the Seiberg-Witten equations 
on $Y \times \R$ with the pull-back $\spinc$ structure, modulo 
gauge transformations. For a general introduction to Seiberg--Witten
Floer theory see \cite{MW}.

Throughout the paper we are considering an oriented, closed homology 
3-sphere $Y$ and a knot $K$ smoothly embedded in $Y$. We consider 
two other manifolds obtained by Dehn surgery on $K$: a homology 3-sphere 
$Y_1$, obtained by $+1$-surgery on $K$, and a 3-manifold $Y_0$ which 
has the homology of $S^1 \times S^2$, obtained by $0$-surgery on $K$. 
Our main goal is to establish the existence of an exact triangle 
relating the Seiberg-Witten-Floer homology groups of these manifolds. 
A similar setup for instanton homology in Donaldson theory was 
considered in \cite{BD}, where Floer's ideas on the corresponding 
construction of the exact triangle for instanton homology are presented. 
 
Because of various technical difficulties intrinsic in this 
program, we need to subdivide the problem into several steps. In 
this first paper we deal with the ``geometric triangle'', namely we 
introduce a suitable ``surgery perturbation'' $\mu$ for the Seiberg-Witten 
equations on $Y$ that simulates the effect of surgery.  
We use the notation ${\cal M}_{Y,\mu}$ for 
 the moduli space of gauge classes of 
solutions of the perturbed Seiberg-Witten equations on $Y$, 
${\cal M}_{Y_1}$ and ${\cal M}_{Y_0}(\s)$ for the moduli spaces
of the perturbed Seiberg-Witten monopoles on
$Y_1$ and $(Y_1, \s)$, where $\s$ is a $\spinc$ structure on $Y_0$. 
 
Our main result  
in this paper is to prove the following  
decomposition theorem for ${\cal M}_{Y,\mu}$. 
 
\begin{The} 
\label{surgery-SW3} 
With a careful choice of perturbations and metrics on $Y, Y_1$ 
and $Y_0$, we have the following relation between the critical sets of the  
Chern-Simons-Dirac functional on the manifolds $Y, Y_1$ 
and $Y_0$: 
\ba 
{\cal M}_{Y,\mu} \cong {\cal M}_{Y_1}\cup \bigcup_{\s_k} {\cal 
M}_{Y_0}(\s_k),  
\na 
where $\s_k$ runs over the  $\spinc$-structures on $Y_0$.  
\end{The} 
 
In section 2, we will briefly review the perturbation theory we use to
define our moduli spaces. In this paper, we only introduce
perturbations sufficient to achieve transversality of moduli spaces
of critical points. Eventually, when dealing with the full
Seiberg-Witten-Floer homology, we shall need a more sophisticated class of
perturbations that achieve transversality simultaneously for moduli
spaces of critical points and of flow lines. These will be non--local
perturbations of the Chern--Simons--Dirac functional, somewhat similar to
those proposed in \cite{K}. We shall deal with this more refined
perturbation theory elsewhere.

In section 3, we will study the 
Seiberg-Witten monopoles on the knot complement $V= Y-K$, 
equipped with a cylindrical end metric modelled on 
$T^2 \times [0, \infty )$.  
We use the notation $\chi(T^2, V)$ 
for the moduli space of  
flat connections on $T^2$ modulo the subgroup of gauge 
transformations on $T^2$ which 
can be extended to $V$. 
Notice that $\chi(T^2, V)$ is a $\Z$-covering  
of the moduli space of flat connections 
on $T^2$ modulo the gauge group $Map(T^2, U(1))$ 
which we denote by $\chi(T^2)$. In 
$ \chi (T^2)$, there is a unique point $\Theta$ such that  
the Dirac operator on $T^2$ coupled with $\Theta$ has non-trivial 
kernel. The main result in section 2 is the following 
structure theorem for the monopole moduli space $\M_V$. 
 
\begin{The}\label{structure} 
For generic metrics and perturbations, 
the moduli space of Seiberg-Witten 
monopoles on $V$, 
denoted by $\M_V$, 
consists of the union of a circle of reducibles $\chi (V)$ and  
an 
irreducible piece $\M_V ^*$ which is 
a smooth oriented 1-dimensional manifold, compact except  
for finitely many ends limiting to $\chi (V)$. 
Moreover, there is a continuous boundary  
value map 
\ba 
\M_V \stackrel{\d_\infty}{\to} \chi(T^2, V) \stackrel{\pi}{\to} \chi 
(T^2). 
\label{map} 
\na 
defined by taking 
the asymptotic limit of the Seiberg-Witten 
monopoles on $V$ over the end. Under $\d_\infty$, $\chi (V)$ is mapped to 
a circle in $ \chi(T^2, V)$, and the compactification $\bar \M_V ^*$
of $\M_V ^*$  is mapped to a collection  
of compact immersed curves in $\chi(T^2, V)$ whose 
boundary points consist of  
a finite set of points 
in $\pi^{-1}(\Theta) \cup  
\d_\infty (\chi (V))$. For  
generic perturbations 
the interior of the curve $\d_\infty (\M_V ^*)$ is transverse to any  
given finite set of curves in 
$ \chi(T^2, V)$. 
\end{The} 

For simplicity of notation, in the following we shall not distinguish
between $\chi(V)$ and its embedded image $\partial_\infty(\chi(V))
\subset \chi(T^2,V)$.
 
In section 4, we will establish a gluing theorem for the moduli  
spaces of critical 
points of the Chern-Simons-Dirac functional when cutting and gluing 
the 3-manifold along a torus. In our case, these are the 
moduli spaces of monopoles on a closed manifold which is  
either $Y$, $Y_1$, or $Y_0$. Let $\nu(K)$ be a tubular neighbourhood 
of $K$ in a closed manifold $Z$, so $Z=V \cup \nu(K)$. 
We may cut $Z$ along $T^2$ and glue in a long cylinder $[-r, r] \times T^2$, 
resulting in a new manifold denoted by $Z(r)$. Use the notation  
$\chi (T^2, Z)$ 
for the character variety (or moduli space) of flat connections on a trivial 
line bundle over $T^2$ modulo the gauge transformations on $T^2$ which can be 
extended to $Z$. 
We denote by $\chi (\nu(K),  Z)$, 
the moduli space of flat connections 
on $\nu(K)$  modulo the gauge transformations on $\nu(K)$ which can be 
extended to $Z$. there is a natural map 
$\chi (\nu(K), Z) \to \chi (T^2, Z)$. We denote by
$\M_{V,Z}^*$ the moduli space of the 
Seiberg-Witten monopoles on $V$ modulo the gauge transformations on $V$ 
which can be extended to $Z$. 
We have a refinement of the  boundary value map of (\ref{map}): 
\ba 
\M^*_{V,Z} \longrightarrow \chi(T^2, Z).   
\label{refine-d-map} 
\na 
Then we have the following gluing theorem. 
 
\begin{The} \label{Gluing} 
 For a sufficiently large $r$, under suitable perturbations and choice of
 metrics,  
there exist the following diffeomorphisms given by the 
 gluing maps on the fibered products 
$$ \#_Y : {\cal M}^*_{V, Y} 
\times_{\chi(T^2,Y)}\chi (\nu(K), Y) 
\longrightarrow {\cal M}^*_{Y(r)}, $$ 
$$ \#_{Y_1} : {\cal M}^*_{V, Y_1}  
\times_{\chi(T^2,Y_1)}\chi (\nu(K), Y_1) 
\longrightarrow {\cal M}^*_{Y_1(r)}, $$ 
$$ \#_{Y_0} : {\cal M}^*_{V, Y_0}  
\times_{\chi(T^2,Y_0)}\chi (\nu(K), Y_0) 
\longrightarrow \bigcup_{\s}{\cal M}^*_{Y_0(r)}(\s). $$ 
Here, ${\cal M}^*_{Y(r)}$, ${\cal M}^*_{Y_1(r)}$ and 
${\cal M}^*_{Y_0(r)}$ are the moduli spaces of irreducible monopoles 
on $Y(r), Y_1(r)$ and $Y_0(r)$ respectively,  
and $\s$ runs over all the possible $\spinc$ structures 
on $Y_0$. The fiber product is taken with respect to the 
refined boundary value maps (\ref{refine-d-map}) from $\M^*_{V, Y}, 
\M^*_{V, Y_1}$ and $\M^*_{V, Y_0}$ 
to  $\chi(T^2,Y)$, $\chi(T^2,Y_1)$ 
and $\chi(T^2,Y_0)$ respectively. 
\end{The}

The proof of Theorem \ref{Gluing} is based on balancing the 
slow decay of certain eigenfunctions of the 
linearization at the approximate solutions, against the 
exponential decay of the finite energy solutions on $V$ with 
non-degenerate asymptotic value, thus obtaining an unobstructed gluing.
 
Using the gluing Theorem \ref{Gluing}, together with the construction of the 
perturbation $\mu$ that ``simulates the effect of surgery'',  
we will be able to derive a corresponding  
deformation of the moduli spaces, and the expected relation 
between the generators of the Floer groups as in Theorem \ref{surgery-SW3}. 
 
In the last section, we apply the result of  
Capell-Lee-Miller on the decomposition of spectral 
flow (Theorem C of \cite{CLM2}) to study the relative gradings of 
monopoles under the identification of Theorem \ref{surgery-SW3}. 
We show that the identification of Theorem \ref{surgery-SW3} 
is compatible with the relative gradings on 
the Seiberg-Witten-Floer chain complexes (Cf. Proposition 
\ref{same-deg-1-mu} and Proposition \ref{same-deg-0-mu}). 
 
{\bf Acknowledgments} 
 We are very grateful to Ronnie Lee, Tom Mrowka, Vicente Munoz, 
 and Peter Ozsvath for useful discussions and
suggestions. We like to thank Cliff Taubes for providing 
the proof of Lemma \ref{taubeslemma} and Liviu
Nicolaescu for the proof of Lemma \ref{Liviu-metrics}.  
The three authors also thank the Max-Planck-Institut f\"ur 
Mathematik, Bonn for the kind hospitality and for support. 
AC and BW are partially supported by Australian Research Council. 
MM is partially supported by NSF grant DMS-9802480 and by Humboldt
Foundation (Sofja Kovalevskaya Award).

\section{Seiberg-Witten equations on 3-manifold}

The 3-dimensional Seiberg-Witten monopoles on a compact manifold
have been extensively studied in \cite{CMWZ} \cite{Chen}
\cite{Froy} \cite{KM1} \cite{KM2} \cite{MW} \cite{MT} \cite{MST}. In this
section we will briefly recall some of the main features of
3-dimensional monopoles. A 3-dimensional monopole, as noted first
in \cite{KM1}, can be viewed as a critical point of the Chern-Simon-Dirac
functional on an infinite dimensional space (the orbit space of
$\spinc$ connections and sections of the spinor bundle under
the action of the gauge group). We recall the basic setting of
3--dimensional Seiberg--Witten theory, then we will end this section
with the observation that, under a generic perturbation with compact
support in a fixed open set, the critical points are all non-degenerate.

Let $(Y, g)$ be a closed, oriented Riemannian 3-manifold. A $\spinc$ structure
$\s$ on $(Y, g)$ is a pair $(W, \rho)$ consisting of a
rank 2 Hermitian bundle $W$ together with a Clifford multiplication
$\rho: T^*Y \to End (W)$. If $\{e_1, e_2, e_3\}$ are an oriented orthonormal
frame for $TY$, we choose the Clifford multiplication such that
$\rho(e_1)\rho(e_2)\rho(e_3) =1$.

With the Levi-Civita connection on the frame bundle of $Y$, 
a $U(1)$-connection $A$ 
on the determinant bundle $det (W)$ determines a $\spinc$ connection
$\nabla_A$ on $W$ such that $\rho$ is parallel. 
Applying  Clifford multiplication,
we can define a Dirac operator, denoted by $\dirac _A$. Then the
Seiberg-Witten equations are the
equations for a pair $(A, \psi)$ consisting of a $U(1)$-connection
on $det (W)$ and a section $\psi$ of $W$ ($\psi$ is called
a spinor):
\ba
\left\{ \begin{array}{l}
*F_A = \sigma (\psi, \psi) + \mu \\[2mm]
\dirac_A (\psi)=0.
\end{array}\right.
\label{3dSW}
\na
Here $\mu$ is a co-closed imaginary-valued 1-form on $Y$, and 
$\sigma (\cdot, \cdot)$ is a symmetric $\R$-bilinear form 
$W\otimes W \to T^*Y \otimes i\R$ given by
\[\begin{array}{lll}
\sigma (\psi, \psi) &=& -\rho^{-1}((\psi\otimes \psi^*)_0)
= -\rho^{-1}(\psi\otimes \psi^* - \displaystyle{\frac {|\psi|^2}{2}} Id)\\[2mm]
&=& \displaystyle{\frac i 2} Im \langle  \rho(e_i)\psi, \psi \rangle  e^i
\end{array}.
\]
Note that this  $\R$-bilinear form $\sigma (\cdot, \cdot)$ satisfies the
following property \cite{CW}:

(1) Under Clifford multiplication, we have  $\sigma (\psi, \psi).\psi
= -  \displaystyle{\frac 12} |\psi|^2 \psi$, and $\langle
\alpha. \psi, \psi\rangle  
= 2 \langle \alpha, \sigma (\psi, \psi)\rangle _{T^*Y}$, for
$\alpha \in \Omega^1(Y, i\R)$.

(2) $\sigma (\psi, \phi) =0$ if and only if on $Y-\psi^{-1}(0)$
 $\phi = i r \psi$
for a real-valued function $r$ on $Y-\psi^{-1}(0)$.

(3) For any imaginary valued 1-form $\alpha$, 
$\sigma (\alpha. \psi, \phi) + \sigma (\psi, \alpha.\phi) 
= - (Re \langle \psi, \phi\rangle ) \alpha$.  

(4) If $\psi$ is a nowhere vanishing section of  $W$, then
$W \cong \underline{ \CC} \psi \oplus \psi^{\perp}$, and
$\sigma (\psi, \cdot)$ defines a bundle
isomorphism between
$ \underline{ \R}  \psi \oplus  \psi^{\perp}$ and
$T^*Y \otimes i\R$.

Denote by ${\cal A}_Y$ the configuration space of $(Y, \s)$ consisting
of pairs $(A, \psi)$ with the completion under $L^2_1$-norm. The
gauge group of automorphisms of the $\spinc$-bundle $W$ is 
$\G_Y = Map(Y, U(1))$ with $L_2^2$-completion. $\G_Y $ acts on
${\cal A}_Y$ by
\[
u(A, \psi) = (A-2u^{-1}du, u\psi),\]
and the Seiberg-Witten equations are invariant under this
action. Denote by $\B_Y$ the quotient space of ${\cal A}_Y$ by the gauge
group action. $\B_Y$  is an infinite dimensional Hilbert manifold
except at points where the spinor part is zero, which are
called reducible points. Otherwise, points $(A, \psi)$ with $\psi\neq 0$
are called irreducible. As noted in \cite{KM1}, the 
Seiberg-Witten equations on $(Y, \s, g)$ are the equations for
the critical points of the following Chern-Simons-Dirac
functional on ${\cal A}_Y$:
\ba
{\cal C}_\mu (A, \psi) = -\displaystyle{\frac {1}{2} \int_Y} (A-A_0) 
\wedge (F_A +F_{A_0} -2*\mu) + 
\displaystyle{\int_Y} \la \psi, \dirac_A \psi\ra dvol_Y,
\label{CSD}
\na
where $A_0$ is a fixed connection on $det(W)$. Note that
$\C_\mu$ descends to a circle-valued function on $\B_Y$. The
set of critical points of $\C_\mu$  on $\B_Y$ is denoted by
$\M_{Y, \mu}(\s)$, its irreducible critical point set
 is denoted by $\M_{Y, \mu}(\s)^*$.

For any critical point
$(A, \psi)$ on ${\cal A}_Y$, the infinitesimal action of $\G_Y$ and the derivative
of $grad (\C_\mu)$ at $(A, \psi)$ define a complex
\ba
\Omega^0_{L^2_2}(Y, i\R) \stackrel{G_{(A, \psi)}}{\to} 
\Omega^1_{L^2_1}(Y, i\R) \oplus L^2_1(W) \stackrel{H_{(A, \psi)}}{\to}
\Omega^1_{L^2}(Y, i\R) \oplus L^2(W),
\label{3d-complex}
\na
where the maps $G_{(A, \psi)}$ and $H_{(A, \psi)}$
are given by 
\[
\begin{array}{l}
G_{(A, \psi)} (f) = (-2df, f\psi), \\[2mm]
H_{(A, \psi)} (\alpha, \phi ) = (*d\alpha - 2 \sigma (\psi, \phi),
\dirac_A \phi + \displaystyle{\frac 12} \alpha. \psi).
\end{array}
\]

We say that $[A, \psi]$ is a non-degenerate critical point
of $\C_\mu$ on $\B_Y$ if the middle cohomology of (\ref{3d-complex})
is zero:
\[
Ker H_{(A, \psi)} / Im G_{(A, \psi)}  =0.
\]
At the smooth points of $\B_Y$, this definition is the same as
saying that the derivative of $grad (\C_\mu)$ at a critical point
is non-degenerate. The gradient of $\C_\mu$
can be viewed as  an $L^2$-tangent vector field on $\B_Y$, a section of
the $L^2$-tangent bundle over $\B_Y$, while
the tangent space of $\B_Y$ at $[A, \psi]$ is the
$L^2_1$-completion of 
\[
Ker G_{(A, \psi)}^* =\{ (\alpha, \phi)| d^* \alpha + iIm \langle  \psi, \phi\rangle  =0.\}
\]

The covariant derivative of $grad (\C_\mu)$, denoted by $H_{[A, \psi]}$, 
defines a operator on $Ker G_{(A, \psi)}^*$, sending
$(\alpha, \phi) \in Ker G_{(A, \psi)}^*$ to
\[
(*d \alpha -2 \sigma (\psi, \phi) -2df, \dirac_A\phi + \frac 12 \alpha.\psi
+ f\psi),
\]
where $f$ is the unique solution to the equation
\[ (d^*d  + \frac 12 |\psi|^2 )f = iIm \langle  \dirac_A\psi, \phi\rangle .\]
Note that $H_{[A, \psi]}$ is a closed, unbounded, essentially self-adjoint, 
Fredholm operator on the
$L^2$-completion of $Ker G^*_{(A, \psi)}$, its eigenvectors form
an $L^2$-complete orthonormal
basis, its $L^2$-spectrum forms a discrete subset of the
real line with no accumulation points. Hence, as in \cite{MW},
the spectral flow of $H_{[A, \psi]}$, along a path connecting
two critical points defines a relative index on 
 $\M^*_{Y, \mu}(\s) \times \M^*_{Y, \mu}(\s) $.
This relative index depends only on the homotopy class of the connecting 
path for non-torsion $\spinc$ structure (Cf. Remark 4.5 in \cite{MW}
and Definition 3.6 in \cite{CW}).

The following properties about the critical points of $\C_\mu$ are
now standard (See \cite{MW} \cite{Chen} \cite{Froy} \cite{KM2}
\cite{Lim2}). 

\begin{Pro} \label{monopole-3d} 
There exists a Baire set of co-closed 1-form $\mu
\in \Omega^1_{L^2_2}(Y, i\R)$ such that all the critical points in 
$\M_{Y, \mu}(\s)$ are non-degenerate. Moreover, if $b_1(Y) >0$, 
$\M_{Y, \mu}(\s)$ consists of only finitely many irreducible
points in $\B_Y$; if $Y$ is a rational homology 3-sphere, assume
that  a generic $\mu$ satisfies $Ker \dirac_\theta =0$
(where $\theta$ is the unique reducible point in $\M_{Y, \mu}(\s)$,
that is, $*F_\theta =\mu$), then
$\M_{Y, \mu}(\s)^* = \M_{Y, \mu}(\s) - \{\theta\}$
 consists of only finitely many irreducible
points.
\end{Pro}

In this paper and sequel work, it is convenient to use a perturbation
with support contained in a fixed open set, so that Proposition
\ref{monopole-3d} still holds for perturbations with 
compact support contained in a fixed open set.
The first such statement was made 
in Proposition 7.1 \cite{Taubes} by Taubes, who
kindly communicated the proof to us.

\begin{Pro} \label{3dim-compactperturbation}
Fix a non-empty open set $U$ in $Y$ and a $\spinc$ structure $\s$ on $Y$,
if $b_1(Y)>0$ and $c_1(det (\s)) =0$, we require that $U$
is chosen so that the map $H^2(Y, \R) \to H^2(U, \R)$ is non-zero. Then 
there exists a Baire
set of co-closed imaginary valued 1-forms $\mu$ with compact support
in $U$ such that all the critical points of ${\cal C}_\mu$ on
$\B_Y$ are non-degenerate.
\end{Pro}
\begin{proof} 
We first study the family version of the critical points
of ${\cal C}_\mu$ on $\B^*_Y$, where $\mu$ is from a set of imaginary valued
co-closed 1-forms on $Y$ with compact support in $U$. Denote this
set of perturbations as $Z(U, i\R)$. Let $[\mu, A, \psi]$
be a critical point of ${\cal C}_\mu$. We need to show that the
derivative of the gradient of $\{{\cal C}_\mu\}_{\mu \in Z(U, i\R)}$ 
is surjective. Namely, consider
\[
KerG^*_{[A, \psi]} \times Z(U, i\R) \to KerG^*_{[A, \psi]},
\]
which sends $ (\alpha, \phi, \mu_1)$ to
\[
(*d\alpha -2 \sigma (\psi, \phi) + \mu_1, \dirac_A\phi + \frac 12
\alpha.\psi).\]

Suppose that $(\alpha,  \phi)$ is orthogonal to the image of
the above map, then $(\alpha,  \phi)$ satisfies the following equations:
\ba\begin{array}{lll}
&(1)&  d^* \alpha + iIm \langle  \psi, \phi\rangle  =0,\\[2mm]
&(2)& *d\alpha -2 \sigma (\psi, \phi) =0,\\[2mm]
&(3)& \dirac_A\phi + \frac 12\alpha.\psi=0, \\[2mm]
&(4)& \text{$\alpha$ is exact when restricted to $U$}.
\end{array}\label{3d-surj}
\na

The elliptic regularity implies that $(\alpha,  \phi)$
is smooth.  From (4) and (2) of equations (\ref{3d-surj}), we know
that $\sigma (\psi, \phi) =0 $ on $U$. The following Lemma due to
Taubes \cite{TauPriv} 
will ensure that $\sigma (\psi, \phi) =0 $ on $Y$. Hence, there
is a real-valued smooth function f on $Y$,
such that $\phi = if\psi$. Using (3) of (\ref{3d-surj}),
we obtain 
\[
\dirac_A (if\psi) + \frac 12\alpha.\psi =0
\]
on $Y$, which leads  to $\alpha = -2i df$ on $Y$. By the equation
(1) in (\ref{3d-surj}), we get
\[
2d^*d f + f |\psi|^2 =0 \]
on Y.  Note that $\psi^{-1}(0)$ does not disconnect any domain in $Y$
(the unique 
continuation principle for Dirac operator (see page 57-58 
\cite{FU})).  Therefore, $f=0 $ which implies that
$(\alpha, \phi) =0$. 

{}From the Sard-Smale theorem, there is a Baire set of $\mu \in Z(U, i\R)$
such that all critical points of ${\cal C}_\mu$ in $\B^*_Y$ are
non--degenerate for a generic $\mu$.

Now we need to prove that the reducible critical point of ${\cal C}_\mu$
is also non-degenerate. By the assumption, ${\cal C}_\mu$ admits
reducible critical point if and only if
$Y$ is a rational homology 3-sphere. From the analysis in
\cite{MW}, we know that, in order to achieve the non-degenerate
condition at reducible critical point,  $\mu$ is required
to be away from the codimension
one subset $Z(U, i\R)$ where the corresponding Dirac operator
has non-trivial kernel. This completes the proof of the Proposition.
Now we give the proof of Taubes' Lemma.  

\begin{Lem} \label{taubeslemma}
 (Taubes) Let $(A, \psi)$ and $(\alpha, \phi)$ as above,
where $(A, \psi)$ is a solution to the Seiberg-Witten equation (\ref{3dSW})
and $(\alpha, \phi)$ satisfies (1)-(3) of the equations (\ref{3d-surj}). Then
$q= \sigma(\psi, \phi)$ obeys an equation of the form
\[
\Delta q = H\cdot q + K \cdot \nabla q
\]
at all points where $\psi\neq 0$. Here $\Delta$ is
the Laplacian on differential 1-forms and $H$ and $K$
are linear maps that depend implicitly on $\psi$. 
The set of points where $\psi\neq 0$ is path connected open
dense set in $Y$. The unique
continuation principle applies to $q$ so that $q$
cannot vanish on $U$ without vanishing everywhere on $Y$.
\end{Lem}

\noindent \underline{Proof of Taubes' Lemma}: Apply the Laplacian to
$q= \sigma(\psi, \phi)$, we have the following expression of $\Delta q$:
\ba
\Delta q = \sigma( \Delta \psi,  \phi )
+  \sigma(  \psi, \Delta  \phi ) + 2 \sigma( \{\nabla_A \psi, \nabla_A\phi \}_
{T^*Y}),
\label{laplace}
\na
here $\Delta$ acting on spinors is $\nabla_A^*\nabla_A$, and
$\{\nabla_A \psi, \nabla_A\phi \}_
{T^*Y}$ is the pairing using the metric on $T^*Y$. Now invoke 
the Weitzenb\"ock formula for the Dirac operator,
\[
\dirac^*_A\dirac_A = - \Delta + \frac \kappa 4 -\frac 12 \rho (*F_A),
\]
where $\kappa$ is the scalar curvature on $Y$. Thus,
from the Dirac equations for $\psi$ and $\phi$, we obtain
\[
\begin{array}{l}
\Delta \psi = \frac \kappa 4 \psi - \frac 12 (*F_A).\psi, \\[2mm]
\Delta \phi = \frac \kappa 4 \phi  - \frac 12 (*F_A).\phi + 
\frac 12 (d^*\alpha)\psi -  \frac 12 (*d\alpha).\psi - \nabla_A^\alpha\psi ,
\end{array}
\]
here $\nabla_A^\alpha\psi = \{\alpha, \nabla_A\psi \}_{T^*Y}$. Plug
these two equations into (\ref{laplace}), and note that
$\sigma(\psi, \frac 12 (d^*\alpha)\psi) =0$ and 
$\sigma(\psi, -\frac 12 (*d\alpha).\psi) = |\psi|^2 q$.
 We get
\ba\begin{array}{lll}
\Delta  q &=& (\frac \kappa 2 + |\psi|^2)q + \sigma (- \frac 12 (*F_A).\psi,
\phi) + \sigma (\psi, - \frac 12 (*F_A).\phi) \\[2mm]
&& + 2\sigma(\{\nabla_A\psi, \nabla_A\phi \}_{T^*Y})+ \sigma (\psi, 
-\nabla_A^\alpha \psi) \\[2mm]
&=& (\frac \kappa 2 + |\psi|^2)q + *F_A (Re \langle \psi, \phi\rangle ) \\[2mm]
&& + 2\sigma(\{\nabla_A\psi, \nabla_A\phi \}_{T^*Y})+ \sigma (\psi, 
-\nabla_A^\alpha \psi).
\end{array}\label{laplace2}
\na
Write $\phi = i r \psi + \lambda$ where $r$ is a real-valued function
on $Y$ and $Re\langle \lambda, i\psi\rangle  =0$, then  
\[
\begin{array}{c}
Re \langle \psi, \phi\rangle  = Re \langle \psi, \lambda \rangle , \\[2mm]
\sigma(\{\nabla_A\psi, \nabla_A\phi \}_{T^*Y}) = 
\sigma(\{\nabla_A\psi, \nabla_A\lambda\}_{T^*Y})
+ \sigma(\{\nabla_A\psi, idr\otimes \psi\}_{T^*Y}).
\end{array}
\]
Hence the equation (\ref{laplace2}) can be written as
\ba\label{laplace3}
\begin{array}{lll}
\Delta  q 
&=& (\frac  \kappa 2 + |\psi|^2)q +*F_A (Re \langle \psi, \lambda\rangle )
\\[2mm]
&&+ 2\sigma(\{\nabla_A\psi, \nabla_A\lambda\}_{T^*Y}) 
+ \sigma(\psi, -\nabla_A^{\alpha+2idr }\psi).
\end{array}
\na

To complete the proof, we only need to show that $\lambda,
\nabla_A\lambda$ and  
$\alpha+2idr$ can be written as combinations of linear maps
on $q$ and $\nabla q$. On the set of points where
$\psi\neq 0$, $\Omega = Y - \psi^{-1}(0)$, we write
$\psi =|\psi|\tau_1$ where $\tau_1$ is a unit-length spinor.
Choose a local basis $\{\tau_1, \tau_2\}$ for the $\spinc$ bundle,
so that  Clifford multiplication in the local orthonormal
coframe $\{e^1, e^2, e^3\}$ for $T^*Y$ is given by 
\[
\rho(e_1) = \left(\begin{array}{lr}
i &0 \\ 0 & -i \end{array}\right),
\rho(e_2)= \left(\begin{array}{lr}
0 &-1\\ 1 &0 \end{array}\right),
\rho(e_3)= \left(\begin{array}{lr}
0 &  i\\  i &0  \end{array}\right).
\]
where $\{e^1, e^2, e^3\}$ can be expressed as
\[
e^1 = -2i \sigma (\tau_1, \tau_1), e^2= 2i \sigma (\tau_1, i\tau_1),
e^3 = -2i \sigma (\tau_1, \tau_2).
\]

Write $\lambda = u\tau_1 + v\tau_2$ for a real-valued
function $u$ and a complex-valued function $v$, then
\[
q= \sigma (\psi, \phi) = \frac i2 |\psi| (u e^1 -Im (v) e^2 + Re (v)e^3).
\]

On $\Omega$, $\sigma_\psi = \sigma(\psi, \cdot)$
defines a bundle isomorphism between
\[
\underline{\R}.\psi \oplus (\underline{\CC}. \psi)^\perp 
\rightarrow T^*Y \otimes i\R.
\]
Thus, we obtain that $\lambda = \sigma^{-1}_\psi (q)$, and
\[
\nabla \lambda = (\nabla (\sigma^{-1}_\psi)) (q) + \sigma^{-1}_\psi
(\nabla q).
\]

Let $b = a + 2idr $, then from (3) of the equations (\ref{3d-surj}), we have
\[
b.\psi = -2 \dirac_A \lambda,\]
as $\lambda$ can be written in terms of $q$ and $\nabla q$, so
is $b$. This completes the proof of Taubes' Lemma. 
\end{proof}

\section{Monopoles on a 3-manifold with a cylindrical end} 
 
In this section we use techniques developed in \cite{MMR} to 
study  the moduli space of Seiberg-Witten monopoles on the knot complement 
$V$ endowed with an infinite cylindrical end $T^2\times [0,\infty)$. 
Our main aim is to present the proof of Theorem \ref{surgery-SW3}. 
Before we give details, we present an overview of the section introducing  
notation. 
 
Consider the three-manifolds $V$ and $\nu (K)$, respectively the knot 
complement and the tubular neighbourhood of the knot $K$ in the 
homology sphere $Y$. Both are 3-manifolds with boundary a torus $T^2$. 
Equip $V$ with a cylindrical end metric 
and a $\spinc$-structure with trivial determinant 
along the half cylinder $T^2\times [0,\infty)$. 
On $T^2$ we use the standard flat metric induced from $\R^2$.  
 
The perturbed Seiberg-Witten equations on $(V, \s)$ are the equations 
\ba 
\left\{ \begin{array}{l} 
*F_A = \sigma (\psi, \psi) +\mu,\\
\dirac_A \psi= 0,
\end{array}\right. 
\label{SW-YT2} 
\na 
for a pair $(A, \psi)$ consisting of a $L^2_{1, loc}$ $U(1)$ connection 
on $det (W)$ and a $L^2_{1, loc}$ spinor section $\psi$ of $W$.  
the perturbation term $\mu$ is a co-closed and imaginary value 1-form
with compact support contained in a fixed open set $U \subset
 V-(T^2 \times [0, \infty))$.  
We denote the corresponding class of perturbations by $Z(U, i\R)$
 
We define  the  energy of  any  Seiberg-Witten monopole $(A, \psi )$ to be 
\ba
\int_V |F_A|^2 d vol_V < \infty,
\label{finite-energy}
\na
Let $\M_V$ denote the Seiberg-Witten moduli space of solutions of the 
equations (\ref{SW-YT2}) with finite energy condition 
 modulo the gauge transformations ${\cal G}_V= Map_{L^2_{2, loc}}(V,U(1))$. 
 
The flat connections on the determinant bundle,  
modulo the even gauge group $\G_{T^2} 
= Map (T^2, U(1))$, form a torus 
$$ \chi (T^2)= H^1(T^2,\R) / 2 H^1(T^2,\Z), $$ 
which is a $\Z_2 \times \Z_2$ cover of the standard torus 
$Hom(\pi_1(T^2),U(1))=\R^2\slash\Z^2$. 
Let $\chi(T^2, V)$ be the moduli space of flat connections 
modulo the subgroup of the gauge 
transformations on $T^2$ which 
can be extended to $V$. 
Let $\pi$ denote the quotient map 
$\pi: \chi(T^2, V) \to \chi (T^2)$, which is a $\Z$-covering 
map.

Suppose we are given a smooth  
solution $(A, \psi)$ of the Seiberg-Witten equations, 
satisfying the finite energy condition (\ref{finite-energy}). 
Then we will see that there is  
a choice of a connection $\tilde A$ in the gauge class of $A$ that approaches 
a flat connection on $T^2$, while the spinor $\psi$  vanishes in the limit 
on the cylindrical end. That is,   
if $s$ is the coordinate on $[0, \infty)$, we will show that 
$\lim_{s\to \infty}  (\tilde A, \psi) = (a_\infty, 0)$ 
in the appropriate topology, 
for each finite energy solution $(A, \psi)$ to the Seiberg-Witten 
equations (\ref{SW-YT2}). 
Thus the asymptotic limit of the 
Seiberg-Witten monopoles  
on the manifold $V$ with a cylindrical end defines a boundary value map 
\ba 
\M_V \stackrel{\d_\infty}{\to} \chi(T^2, V) \stackrel{\pi}{\to} \chi 
(T^2).  
\label{boundary-T^2} 
\na 
We will show that, in a suitable topology, 
this boundary value map is well-defined and 
continuous. Then, we will describe the structure of the 
moduli space $\M_V$. 
 
\subsection{Monopoles on $T^2\times [0, \infty)$}

We begin with the investigation of the behaviour of the 
solutions of the Seiberg-Witten equations on the cylindrical end 
$T^2 \times [0, \infty )$. 
Fix a flat background connection $A_0$ on the determinant bundle 
$det (W)$ with asymptotic limit $a_0$.   
 
\begin{Lem}  Choose the coordinate $s\in [0, \infty)$ on the 
cylindrical end $T^2 \times [0, \infty)$. Choose  
the $\spinc$ structure over $T^2\times [0, \infty)$ 
to be the pull-back of the $\spinc$ structure on $T^2$ with trivial 
determinant, induced by the complex structure. 
We can write $(A, \psi)$ as 
\[ 
\left\{ 
\begin{array}{l} 
A=A_0+ a(s) + h(s) d s,\\[2mm] 
\psi = (\alpha (s), \beta (s)) \in \Lambda^{0,0} \oplus \Lambda^{0,1} 
= \Gamma (W). 
\end{array}\right. 
\] 
where $a(s) = a^{1,0} (s) + a^{0,1}(s) \in \Lambda^1(T^2, i\R)$, 
$h (s) \in \Lambda ^0(T^2, i\R)$. Then the Seiberg-Witten equations  
(\ref{SW-YT2}) 
can be written in the form 
\ba 
\left\{\begin{array}{l} 
F_{A_0+ a} =  \displaystyle{\frac{i}{2}} (|\alpha|^2 - |\beta |^2 ) 
\omega,\\[2mm] 
\displaystyle{\frac{\d a^{0,1}(s) }{\d s} }= 
i\bar\alpha \beta + \bar \d h, \\[2mm] 
 \left(\begin{array}{ccc} 
    i(\d_s +h ) && \bar\d^*_{a(s)}\\[1mm] 
   \bar\d_{a(s)} && -i(\d_s +h )\end{array}\right) 
\left(\begin{array}{c}\alpha \\[1mm]   \beta \end{array}\right) 
= 0. 
\end{array}\right. 
\label{SW-T^2R} 
\na 
where $\omega$ is the area 2-form on $T^2$ with 
$\displaystyle{\int_{T^2} } \omega = 1$. 
\end{Lem} 
\begin{proof} 
We may choose a trivialization of the cotangent bundle to 
$T^2\times [0, \infty)$ 
so that, using a full-stop to denote Clifford multiplication by a one 
form, 
we can make the identifications: 
\ba 
d s. = \left(\begin{array}{lr} 
i &0 \\ 0 & -i \end{array}\right), \quad 
dx. = \left(\begin{array}{lr} 
0 &-1\\ 1 &0 \end{array}\right), 
\quad dy. = \left(\begin{array}{lr} 
0 &  i\\  i &0  \end{array}\right), 
\label{cliff} 
\na 
Letting the Hodge * on forms on $T^2\times [0, \infty)$ be denoted by 
$*_3$, then under the 
preceding identifications we have 
\[\begin{array}{c} 
*_3 (\sigma (\psi, \psi)) = \displaystyle{\frac{i}{2}} 
 (|\alpha|^2 - |\beta |^2 ) \omega -i (\alpha \bar  \beta + \bar \alpha 
\beta) 
\wedge d s 
\\[2mm] 
 F_A = F_{A_0+a} + (d_{T^2} h -\d_s a )\wedge d s, \end{array} 
\]  
hence we get 
\[ 
\left\{\begin{array}{l} 
F_{A_0 +a} =  \displaystyle{\frac{i}{2}} (|\alpha|^2 - |\beta |^2 ) 
\omega,\\[2mm] 
\displaystyle{\frac{\d a^{0,1}(s) }{\d s}} = i\bar \alpha \beta + \bar\d h. 
\end{array}\right. 
\] 
The form $a(s) \in \Lambda^1(T^2, i\R)$ is uniquely 
determined,  as an $i\R$-valued 1-form, 
by its $(0,1)$-part $a^{0,1} \in \Lambda^{0,1} (T^2)$. 
Similarly, the Dirac operator on $T^2 \times [0, \infty)$ can be 
expressed as 
\[ 
\dirac_{a(s) +h(s)ds} = \left(\begin{array}{ccc} 
i(\d_s +h ) && \bar\d^*_{a(s)}\\[1mm] 
   \bar\d_{a(s)} && -i(\d_s +h )\end{array}\right). 
\] 
Thus gives the Dirac equation as in the Lemma. 
\end{proof} 
 
Let $(A,\psi)$ be an irreducible solution of the Seiberg-Witten 
equations on the manifold $V$. 
Along the cylindrical end $T^2\times 
[0,\infty)$ we can use Lemma 3.1 
to write the Seiberg-Witten equations in the form 
\[ 
\left\{\begin{array}{l} 
\d_s a^{0,1} = i\bar\alpha \beta + \bar \partial h,\\[2mm] 
\d_s \alpha = i\bar\d^*_{a(s)} \beta - h \alpha,\\[2mm] 
\d_s \beta = -i \bar\d_{a(s)} \alpha - h \beta. 
\end{array}\right. 
\] 
with the constraint  
$ 
F_a =  \displaystyle{\frac{i}{2}} (|\alpha|^2 - |\beta |^2 
) \omega.$ 
These equations are gauge-equivalent, through a gauge transformation
in ${\cal G}_V$,  to the following equations: 
\ba 
\left\{\begin{array}{l} 
\d_s a^{0,1} = i\bar\alpha \beta ,\\[2mm] 
\d_s \alpha = i\bar\d^*_{a(s)} \beta,\\[2mm] 
\d_s \beta = -i \bar\d_{a(s)} \alpha, 
\end{array}\right. 
\label{SW-T2} 
\na 
on the configuration space ${\cal A}_{T^2}$ 
of triples $(a,\alpha,\beta)$, where $a$ is a $U(1)$-connection 
on $det (W)$ and $(\alpha,\beta)$ is a section 
of $\spinc$ bundle $W$ over $T^2$.  The following Lemma shows 
that (\ref{SW-T2}) can be interpreted as a gradient flow equation. 
 
\begin{Lem} The equations (\ref{SW-T2}) are the downward gradient 
flow equations of the $\G_{T^2}= Map (T^2, U(1))$-invariant  functional 
\ba 
f(a, \alpha, \beta) = - \int_{T^2} \la \alpha, i\bar\d^*_a \beta\ra 
\omega \label{funct-f} 
\na 
on the space ${\cal A}_{T^2} $, where the product $\la, \ra$ 
denotes the natural inner product using the Hodge star operator. 
\label{Morse-T2} 
\end{Lem} 
\begin{proof} 
Direct calculation shows that we have 
\[ 
\nabla \ f (a, \alpha, \beta) = (-i\bar\alpha \beta-i \alpha\bar\beta,  
- i\bar\d^*_{a} \beta, i \bar\d_{a} \alpha ). 
\] 
\end{proof}

Critical points of the functional (\ref{funct-f}) 
with the condition $F_a = \displaystyle{\frac{i}{2}} (|\alpha|^2 - |\beta |^2 
) \omega$ are all the elements  
$(a_\infty, 0,0)$, with $a_\infty$ a flat connection. 
This critical point set is denoted by $\chi (T^2)$, the
quotient space of the flat connection by the
even gauge transformation.
If $(A(s), \psi(s))$ is a solution to the
Seiberg-Witten equation on $[0, \infty) \times T^2$
in temporal gauge, then 
\[
(A(s), \psi(s))= (A_0 + a^{1,0}+ a^{0,1},
(\alpha, \beta))
\]
satisfies the gradient  flow equation of $f$ as given by
(\ref{SW-T2}). 

The next few lemmata describe some fundamental properties
about the solution to the
Seiberg-Witten equation on the cylinder over $T^2$
in temporal gauge.

\begin{Lem}\label{variation1}
Let $\gamma (s)= (A(s), \psi(s))$ be  a solution to the
Seiberg-Witten equation on $[s_1, s_2] \times T^2$
in temporal gauge,  then
\[
\int_{s_1}^{s_2} \|\nabla f (\gamma (s)) \|^2 _{L^2(T^2)} d s
= \int_{[s_1, s_2] \times T^2} (|\nab_A \psi |^2 + |F_A|^2) d vol.
\]
\end{Lem}
\begin{proof}
Since $ (A(s), \psi(s))$ satisfies the Seiberg-Witten
equation on $[s_1, s_2] \times T^2$:
\[
\left\{ \begin{array}{l}
\dirac_A \psi = 0, \\
*F_A = \sigma (\psi, \psi)
\end{array}\right.
\]
The Weitzenb\"ock formula for the
Dirac operator on $[s_1, s_2] \times T^2$ with flat metric then gives 
\[
\dirac_A^* \dirac_A \psi = \nab_A^*\nab_A \psi - \frac{1}{2}
(*F_A). \psi =0.\]
Take  the inner product of both sides  with $\psi$, use
the Seiberg-Witten equation, and note that
$\la (-*F_A). \psi, \psi \ra =  2\la *F_A, \sigma (\psi, \psi)\ra =2
|F_A|^2$.  We obtain 
\[
\frac{1}{2} d^*d |\psi|^2 + |\nab_A \psi|^2 + |F_A|^2 =0.
\]
Integrating the above identity over $[s_1, s_2] \times T^2$, we can write
the result as
\[\begin{array}{lll}
&&\displaystyle{\int_{[s_1, s_2] \times T^2} }
 (|\nab_A \psi|^2 + |F_A|^2) d vol\\[3mm]
&=&\displaystyle{ -\frac 12 \int_{[s_1, s_2] \times T^2} }
d*d |\psi|^2 \\[3mm]
&=&\displaystyle{ -\frac 12 \int_{ \partial ([s_1, s_2] \times T^2)}}
(\partial_s \la \psi, \psi \ra) \omega\\[3mm]
&=& \displaystyle{\int_{T^2}} 
\la \alpha (s_1), i\bar\partial_{A(s_1)}^* \beta (s_1)\ra\omega
-\displaystyle{ \int_{T^2}} \la \alpha (s_2), 
i\bar\partial_{A(s_2)}^* \beta (s_2)\ra\omega.
\end{array}
\]
Here we write $\psi (s) = (\alpha (s), \beta(s))$ 
as a spinor on $T^2$ and use the
equation (\ref{SW-T2}) for $\partial_s\psi$. 
Note that $\gamma (s) = (A(s), \alpha (s), \beta(s))$
solves the gradient flow equation of $f$, hence
\[\begin{array}{lll}
&&
\displaystyle{\int_{s_1}^{s_2} }\|\nab f(\gamma (s))\|^2_{L^2(T^2)} ds\\[3mm]
&=&\displaystyle{\int_{T^2}}
 \la \alpha (s_1), i\bar\partial_{A(s_1)}^* \beta (s_1)\ra\omega
- \displaystyle{\int_{T^2}}
\la \alpha (s_2), i\bar\partial_{A(s_2)}^* \beta (s_2)\ra\omega \\[3mm]
&=&
\displaystyle{\int_{[s_1, s_2] \times T^2} }
 (|\nab_A \psi|^2 + |F_A|^2) d vol.
\end{array}
\]
\end{proof}

\begin{Lem}\label{variation2}
Let $\gamma (s)= (A(s), \psi(s))$ be  a solution to the
Seiberg-Witten equation on $N= [t-1, t+1] \times T^2$
in temporal gauge for any $t \in [0, \infty)$. If 
$\displaystyle{\int_N |F_A|^2 dvol  =  E_N }$, then 
there exists a constant $C_0$ such that the following estimates
hold
\[
\displaystyle{\int_{[t-\frac 12, t+\frac 12] \times T^2}} 
| \nab_{A}\psi |^2 d vol \le C_0 \sqrt {E_N};
\]
\[
\displaystyle{\int_{t-\frac 12}^{t+\frac 12} }
\|\nab f(\gamma (s)) \|^2 _{L^2} ds \le C_0 \sqrt {E_N} + E_N.
\]
Moreover, if $(A(s), \psi(s))$ is a solution to the
Seiberg-Witten equation on $[-1, \infty)  \times T^2$
in temporal gauge with finite energy, then
the corresponding flowline on ${\cal A}_{T^2}$ of $f$
has finite variation of $f$ along $[0, \infty)\times T^2$.
\end{Lem}
\begin{proof}
{}From the $L^2$-bound on $|F_A|$, we immediately obtain
a $L^4$-bound on $\psi$ from the Seiberg-Witten equation,
\[
\|\psi\|^4_{L^4(N)}=
\displaystyle{\int_N }|\psi |^4 dvol =\displaystyle{\frac 14} 
\displaystyle{\int_N }|F_A|^2 dvol.
\]
By the Cauchy-Schwartz inequality, we get
\[
\|\psi\|^2_{L^2(N)} \leq \sqrt{Vol(N)} \|\psi\|^2_{L^4(N)}
= \displaystyle{\frac{\sqrt {2}}{2}} \|F_A\|_{L^2(N)}.
\]
Here we use that $Vol (N)=2$. In the proof the previous lemma,
we found that $\psi$ satisfies 
\[
\frac{1}{2} d^*d |\psi|^2 + |\nab_A \psi|^2 + |F_A|^2 =0.
\]
Multiplying both sides of the above equation with
a cut-off function $\rho$ which
equals 1 on $[t-\frac 12, t+\frac 12]\times T^2$ and vanishes
near the boundary of $N$, and then integrating by parts,
we obtain 
\[
\|\nab_A\psi\|^2_{L^2([t-\frac 12, t+\frac 12]\times T^2)}
\leq \displaystyle{-\frac {1}{2} \int_N}d^*d |\psi|^2 \rho dvol 
\leq C_1 \|\psi\|^2_{L^2(N)} \leq \displaystyle{\frac{\sqrt {2}}{2}} C_1
\|F_A\|_{L^2(N)},
\]
where $C_1$ is a constant depending only on the
cut-off function $\rho$. Putting the above inequalities
together we get the estimates as claimed with 
$C_0 =  \displaystyle{\frac{\sqrt {2}}{2}} C_1$.

The finite variation of $f$ along $[0, \infty)\times T^2$
for a solution on $[-1,  \infty)\times T^2$
is the direct consequence of adding up over a sequence of
middle tubes of length $2$, 
namely, $\{[i-1, i+1]\times T^2| i=0, 1, 2, \cdots \}$, hence
\[
\begin{array}{lll}
&& \displaystyle{\int_0^\infty} \|\nab f(\gamma (s))\|^2_{L^2(T^2)}ds\\[3mm]
&\leq & 2  \displaystyle{\int_{-1}^\infty} \|F_A\|^2_{L^2(T^2)}ds
+ 2C_0 \sqrt{\displaystyle{\int_{-1}^\infty} \|F_A\|^2_{L^2(T^2)}ds}
< \infty.
\end{array}
\]  
\end{proof}

\begin{Lem} Let $(a(s), \alpha(s), \beta(s))$ be a monopole
on $T^2\times [0, \infty)$ with finite energy
$$\int_{T^2\times[0, \infty)}
|F| ^2 dvol < \infty. $$
Then there exists a sequence
$\{s_n \}$ such that $lim_{s_n\to \infty} (a(s_n), \alpha(s_n), \beta(s_n))$
exists and represents a point in $\chi (T^2)$.
\end{Lem}
 \begin{proof} Write the curvature of $a(s)$ on $T^2\times [0, \infty)$
as $F_{a(s)} -\partial_s (a(s)) \wedge ds$, where
$F_{a(s)}$ is the curvature on $T^2$. Then we have  the following
calculation:
\[
\begin{array}{lll}
&& \displaystyle{
\int_{T^2\times [0, \infty)} |F| ^2 dvol} \\[3mm]
&=& \displaystyle{
\int_{T^2\times [0, \infty)} (|F_{a(s)}|^2 + |\partial_s (a(s))|^2) dvol},
\end{array}
\]
which implies that as $s\to \infty$,
\[
|F_{a(s)}| \to 0, \qquad |\partial_s (a(s))| \to 0.
\]
By Uhlenbeck's weak compactness result and the compactness of $\chi (T^2)$,
we know that $a(s)$ weakly converges to a flat connection.

By the monopole equation on $T^2\times [0, \infty)$,
we also obtain
\[
|\alpha (s) |^2 - |\beta(s)|^2 \to 0, \qquad \bar\alpha(s) \beta(s)\to 0.
\]
This implies that $[a(s), \alpha (s), \beta(s)]$ converges weakly
to a point in $\chi (T^2)$. 
\end{proof}

To establish  strong convergence to a point in $\chi (T^2)$
for any finite energy monopole
on $[0, \infty) \times T^2$, we need to apply L. Simon's type
result of "small energy implying small length" as in \cite{MMR}.
We will address this issue at the end of this subsection.

Let $\G_{T^2}^0$ be the based gauge group on $det (W)$, that is, those 
gauge transformations which equal the identity at a fixed based point. 
Denote by $\B_{T^2}^0$ the quotient space of ${\cal A}_{T^2}$ by the 
free action of $\G_{T^2}^0$. Note that the gradient flow 
of $f$ preserves the constraint $F_a =
\frac{i}{2}(|\alpha|^2-|\beta^2|)\omega$, hence we can consider
gradient flow lines of $f$ restricted to 
$${\cal C}_{T^2} = \Big\{(a, \alpha, \beta) 
| F_a = \displaystyle{\frac i2}(|\alpha|^2-|\beta^2|)\omega 
\Bigl\}/\G_{T^2}^0$$ 
as a subset of  $\B_{T^2}^0$.  
 
The space ${\cal C}_{T^2}$ is a singular space,  
the singular set consisting of  
$[a, \alpha, \beta]$ where $a$ is a flat connection and 
$(\alpha, \beta )$ is a spinor section satisfying 
the pointwise condition $|\alpha| = |\beta|$. 
We want to study the asymptotic behavior of the finite 
energy  monopole on $T^2 \times [0, \infty )$, that is, the 
asymptotic behavior of the gradient flow of $f$ restricted  
to ${\cal C}_{T^2}$. 

If we consider a neighborhood of ${\cal C}_{T^2}$
in the whole configuration space $\B^0_{T^2}$,
this introduces new critical points 
which consist of the $[\Theta, \alpha, \beta]$, 
with $\Theta \in \chi (T^2)$ and $(\alpha, \beta)$ satisfying 
\[ 
\bar \alpha \beta =\bar \partial_\Theta (\alpha) =\bar
\partial_\Theta^*(\beta)=0.  
\] 
Note that there is a unique point $\Theta \in \chi (T^2)$ with 
$ker (\bar\d_\Theta + \bar\d_\Theta^*)$ non-trivial, which satisfies  
\[ 
        Ker \bar\d_{\Theta}\cong Ker \bar\d_{\Theta}^* \cong \CC. 
\] 
 
Since we are only interested in the behavior of the monopoles 
on $T^2 \times [0, \infty)$, among flowlines of $f$ on  
$\B^0_{T^2}$,  we only study those that flow  
to the critical manifold $\chi (T^2)$. 
The Hessian operator of $f$ at the critical point $[a_\infty, 0, 0]$ 
in $\chi (T^2)$
is given by 
\ba 
Q_{[a_\infty, 0, 0]} (a_1, \alpha_1, \beta_1) = \Bigl(0, -i 
\bar\d^*_{a_\infty} 
\beta_1, i \bar\d_{a_\infty}\alpha_1 \Bigr). 
\label{Hessian-f} 
\na 
where $(a_1, \alpha_1, \beta_1)$ is a $L^2_1$-tangent vector
of $\B^0_{T^2}$ at $[a_\infty, 0, 0]$, that is, $(a_1, \alpha_1, \beta_1)$ 
satisfies the condition $d^*a_1 =0$, and we view
$(0, -i \bar\d^*_{a_\infty} \beta_1, i \bar\d_{a_\infty}\alpha_1 )$
as a $L^2$-tangent vector of $\B^0_{T^2}$ at $[a_\infty, 0, 0]$.
Then the following lemma is obtained by a direct calculation.

\begin{Lem}  
For $a_\infty \ne \Theta$ in $\chi (T^2)$, 
$f$ is non-degenerate  at $a_\infty$ 
in the sense of Morse-Bott, that is, the Hessian 
operator $Q$ at $[a_\infty, 0, 0]$ is non-degenerate in the 
normal direction to the critical manifold in the tangent space of 
$\B^0_{T^2}$
at $[a_\infty, 0, 0]$. At the point $\Theta$, the kernel of  
the Hessian operator is given by  
\[ 
   H^1(T^2, i\R)\oplus  Ker \bar\d_{\Theta}\oplus  
  Ker \bar\d_{\Theta}^* \cong \CC^3. 
\] 
\end{Lem}

Let $U_\Theta$ be a small open neighbourhood of $\Theta$ in $\chi 
(T^2)$. 
For any point $a_\infty \in \chi (T^2)\backslash U_\Theta$,
the spectrum of $Q_{a_\infty}= Q_{[a_\infty, 0, 0]}$ 
 (as a first order elliptic operator (\ref{Hessian-f}))
is discrete, real and without accumulation points.  Let 
$\mu_{a_\infty} >0$ 
be the smallest  absolute value of the non-zero eigenvalues of the Hessian 
operator $Q_{a_\infty}$. 
Now we can establish the decay estimate for the Seiberg-Witten 
monopoles along the cylindrical end of $V$. The first
exponential decay estimate is for a solution to the
Seiberg-Witten equation on $[0, R] \times T^2 $ ($R>1$) which
is near a critical point in $\chi (T^2)$.

\begin{Lem}\label{decay-short}
Suppose that $x(s)= [a(s), \psi(s)]$ is a 
flow line of $f$, corresponding to an irreducible
finite energy monopole on $T^2 \times [0, R]$ 
in temporal gauge. There is a representative
$(A(s), \psi(s))$ which
is gauge equivalent to
 $(a_\infty, 0) + (b, \eta)$, where $[a_\infty, 0] \neq \Theta$. 
There exist positive
constants $\epsilon, \delta, C_1$ such that, if $(b, \eta)$ has 
$L^2_1$ -norm less than $\epsilon$ on any $s$-slice, then
\[
\| (b(s), \eta(s)) \|_{L^2_1(T^2)}  \leq C_1 
\bigl(exp(-\delta s) + exp(-\delta(R-s))\bigr) \]
on any constant $s$-slice ($s\in [0, R]$).
\end{Lem}
\begin{proof}
Write $\lambda = (b, \eta)$, then $\lambda$ satisfies
the following equation:
\[
\partial_s \lambda = Q_{a_\infty} \lambda + n (\lambda),
\]
Here $n(\lambda)$ is  second order in $\lambda$ with
$\|n (\lambda) \|_{L^2(T^2)} \le \epsilon \|\lambda \|_{L^2(T^2)}$ and
$Q_{a_\infty} = Q_{[a_\infty, 0, 0]}$. Note that
the flowline  of $f$ on $\chi (T^2)$ is static, hence we can establish the  
analogous result as Lemma 5.4.1 in \cite{MMR} as follows. 

Let $\lambda_{\pm}$ denote the projection of
$\lambda$ onto the eigenspaces of
$Q_{a_\infty} $ with positive and negative eigenvalues.
Let $\|\lambda_{\pm}\|$ be the functions on $[0, R]$
given by the $L^2(T^2)$-norm on the $s$-slice of $[0, R]\times T^2$.
Then we have
\[
\partial_s \|\lambda_+ \| - (\mu_{a_\infty} -\epsilon) \|\lambda_+ \|
+ \epsilon \|\lambda_-\| \geq 0;
\]
\[
\partial_s \|\lambda_- \|+ (\mu_{a_\infty} -\epsilon) \|\lambda_-\|
- \epsilon \|\lambda_+\| \leq 0.
\]
When $\epsilon << \mu_{a_\infty}$, from the above inequalities
together with the comparison principle (Cf. Lemma 9.4 \cite{Taubes}), 
we obtain that the
$L^2$-norm of $\lambda$ on the $s$-slice 
is decaying exponentially with decay rate
$\delta \le \mu_{a_\infty}/2$. Then the claim of the lemma
follows from the standard bootstrapping argument.
\end{proof}

\begin{Pro} Suppose that $\gamma (s)= [a(s), \psi(s)]$ is an irreducible 
flow line of $f$, corresponding to an irreducible  
finite energy monopole on $T^2 \times [0, \infty)$, 
with asymptotic limit $[a_\infty, 0, 0]$ where $[a_\infty] \ne \Theta 
\in \chi (T^2)$. 
Then, there exist 
gauge representatives $(a(s), \psi(s))$ 
for $\gamma (s)$ and $a_\infty$ for $[a_\infty, 0, 0]$ such that 
$(a(s) -a_\infty, \psi(s))$ decays exponentially along with its 
first derivative as $s\to \infty$. 
\label{decay-f} 
\end{Pro} 
\begin{proof} From Lemma \ref{variation2}, we know
that the variation of $f$ is finite, that is,
\[
\displaystyle{\int_{1}^\infty} \|\nab f (\gamma (s))\|^2_{L^2(T^2)} ds 
\]
is finite. Then  we have the following estimate, 
whose proof is analogous 
to the proof of Lemma 6.14 of \cite{MST}. We sketch the proof here. 
 
\noindent{\bf Claim}: 
There exist constants $E_0$ and $C$ such that for any $R>1$,
and for  $\gamma(s)= (A(s), \psi(s) )$ any solution to
the Seiberg-Witten equation in  temporal gauge on
$[0, R+1]\times T^2$ satisfies
\[
\displaystyle{\int_{0}^{R+1}}  \|\nab f (\gamma (s))\|^2_{L^2(T^2)} ds \leq E_0,
\]
then we have the estimate 
\[ 
\displaystyle{\int_{1}^R}\|\nab f (\gamma (s))\|^2_{L^2_1(T^2)} ds 
 \leq C \displaystyle{\int_{0}^{R+1}}  \|\nab f (\gamma (s))\|^2_{L^2(T^2)} ds. 
\] 
\noindent{\bf Proof of the Claim} 
Let $\gamma (s)= ((A(s), \psi(s))$ be  a solution to the
Seiberg-Witten equation on $ N= [s_1, s_2] \times T^2$
in temporal gauge,  then from Lemma \ref{variation1}, we have
\[
\int_{s_1}^{s_2} \|\nabla f (\gamma (s)) \|^2 _{L^2(T^2)} d s
= \int_{[s_1, s_2] \times T^2} (|\nab_A \psi |^2 + |F_A|^2) d vol.
\]
Denote by 
\[\begin{array}{lll}
E&=& \displaystyle{\int_{s_1}^{s_2} }\|\nabla f (\gamma (s) \|^2 _{L^2(T^2)} d s
\\[3mm]
&=&\displaystyle{ \int_{N}} (|\partial_s A|^2 + |\partial_s \psi|^2) dvol.
\end{array}
\]
Then we have the following estimates
\[
\|F_A\|_{L^2(N)}\leq \sqrt {E}; \qquad \|\psi\|^2_{L^4(N)}=2 \|F_A\|_{L^2(N)}
 \leq 2 \sqrt {E}.
\]

We proceed as in Lemma 6.14 of \cite{MST} and differentiate the
Seiberg-Witten equations to get
\[ d{\d_s A}=*\sigma(\d_s\psi,\psi) \]
\[ \dirac_A \d_s\psi+(\d_s A)\cdot\psi =0. \]
The gauge fixing condition implies that 
\[
d^*( \d_s A) +i Im \la \d_s\psi,\psi \ra =0.
\]
Introduce a cutoff function $\rho$ identically equal to 1 on the
middle third piece of $N$ and vanishes near the boundary such that
$|d\rho|$ is at most $\frac{M}{s_2-s_1}$ where $M$ is a universal
constant. Set $(V, \lambda ) = 
(\rho \d_s A, \rho \d_s \psi)$.
Then we can estimate the quantity 
\[ SW(V, \lambda ) = \bigl(dV -*\sigma(\lambda ,\psi), \dirac_A (\lambda )
+V\cdot\psi, d^*V +i Im \la \lambda ,\psi\ra\bigr)
\]
by
\[ \| SW(V, \lambda )\|^2_{L^2(N) }\leq \frac{C}{(s_2-s_1)^2}
\left( \| \d_s A \|^2_{L^2(N)}+\| \d_s\psi  \|^2_{L^2(N)}\right). \]
Here $C$ is a universal constant depending only on $\rho$.

On the other hand, we can estimate
\[
\begin{array}{lll} 
\| SW(V, \lambda )\|^2_{L^2(N) }&\geq & \frac{1}{2} \left(
\|dV\|^2_{L^2(N) } + \|\dirac_A (\lambda)\|^2_{L^2(N) }
+ \| d^*V\|^2_{L^2(N)}\right)\\[2mm]
&&-2 \left( \|\sigma(\lambda,\psi)\|^2_{L^2(N) }+ \| V\cdot \psi \|^2
_{L^2(N) }  
 \|Im \la \lambda,\psi\ra \|^2_{L^2(N) }\right).
\end{array} \]
Assume that $l= s_2-s_1 \le 1$, then the Sobolev multiplication
theorem and Sobolev embedding theorem imply that
there are constants $C_0$ and $C_1$ such that
\[
\begin{array}{lll}
\|\sigma(\lambda,\psi)\|^2_{L^2(N) }&\leq& 
C_0 \|\lambda\|^2_{L^4(N) }
\|\psi\|^2_{L^4(N) } \\[2mm]
&\leq& C_1 \sqrt{E} \|\lambda\|^2_{L^2_1(N) }.
\end{array}
\]
Similarly, by choosing $C_1$ appropriately, we have
\[
 \| V\cdot \psi \|^2 _{L^2(N) } \leq C_1 \sqrt{E} \|V\|^2_{L^2_1(N) };
\]
\[
\|Im \la \lambda,\psi\ra \|^2_{L^2(N)}\leq C_1
\sqrt{E}\|\lambda\|^2_{L^2_1(N)}. 
\]
These inequalities imply that
\[ 
\begin{array}{lll} &&\|dV\|^2_{L^2(N) } + \|\dirac_A (\lambda)\|^2_{L^2(N) }
+ \| d^*V\|^2_{L^2(N)}\\[2mm]
&\leq& 8 C_1\sqrt{E}\left( \|V\|^2_{L^2_1(N)} + \|\lambda\|^2_{L^2_1(N)}\right)
+ 2 \| SW(V, \lambda )\|^2_{L^2(N) }
\end{array}
\]
Standard estimate for the elliptic operator
$(d+d^*, \dirac_A)$ can be employed to show that
there is a constant $C_2$ such that 
\[\begin{array}{lll}
\|V\|^2_{L^2_1(N)} + \|\lambda\|^2_{L^2_1(N)}&\leq& C_2 
\left( \|dV\|^2_{L^2(N) } + \|\dirac_A (\lambda)\|^2_{L^2(N) }
+ \| d^*V\|^2_{L^2(N) }\right)\\[2mm]
&& + C_2 \left( \|V\|^2_{L^2(N)}+ \|\lambda\|^2_{L^2(N)}\right).
\end{array}\]
The Cauchy-Schwartz inequality and the Sobolev embedding theorem
imply that there exists a constant $C_3$ such that
\[
\left( \|V\|^2_{L^2(N)}+ \|\lambda\|^2_{L^2(N)}\right)
\leq C_3\sqrt{s_2-s_1}
 \left( \|V\|^2_{L^2_1(N)}+ \|\lambda\|^2_{L^2_1(N)}\right).
\]
Put all these inequalities together, we have
\[\begin{array}{lll}
\|V\|^2_{L^2_1(N)} + \|\lambda\|^2_{L^2_1(N)}&\leq&
\displaystyle{\frac {2CC_2}{(s_2-s_1)^2}}
\left( \| \d_s A \|^2_{L^2(N)}+\d_s\psi  \|^2_{L^2(N)}\right)\\[2mm]
&& + 8 C_1C_2\sqrt{E}\left( \|V\|^2_{L^2_1(N)} + \|\lambda\|^2_{L^2_1(N)}\right)
\\[2mm]
&&+ C_2C_3 \sqrt{s_2-s_1}
 \left( \|V\|^2_{L^2_1(N)}+ \|\lambda\|^2_{L^2_1(N)}\right)
\end{array}
\]
Then there is a constant $E_0$ and  a constant $l_0$ satisfying
\[
1-8C_1C_2\sqrt{E_0} -C_2C_3\sqrt{l_0} \geq \displaystyle{\frac 12}
\]
such that if $E\leq E_0$ and $s_2-s_1= l_0$,   there is
a constant $C_4$ with the following estimate
\[
\|V\|^2_{L^2_1(N)} + \|\lambda\|^2_{L^2_1(N)}\leq
\frac{C_4} {l_0^2} \left( \| \d_s A \|^2_{L^2(N)}+ \|
\d_s\psi  \|^2_{L^2(N)}\right). 
\]
Since on the middle third piece $N'$, $V|_{N'} = \d_s A |_{N'}$
and $\lambda|_{N'} = \d_s\psi |_{N'}$,
this implies that for any tube $N=[s_1, s_2] \times T^2$
of length $l_0$ and any solution $(A(s), \psi(s))$ on $N$ of
energy at most $E_0$, we have
\[
\| \d_s A \|^2_{L^2_1(N')}+ \|\d_s\psi\|^2_{L^2_1(N')} 
\leq \frac{C_4} {l_0^2} \left( \| \d_s A \|^2_{L^2(N)}+ \|
\d_s\psi  \|^2_{L^2(N)}\right). 
\]
Then the estimate in the claim follows by adding up a
sequence of middle third pieces of tubes (length $l_0$) with the
constant $C= 3C_4/l_0^2$ and $E_0$ as above.

With this claim and Lemma \ref{decay-short}, we can prove Proposition
\ref{decay-f} using the method of the proof of Proposition 6.16
\cite{MST} and the fact that $f$ is a Morse-Bott function on $\B^0_{T^2}$
and satisfies the Palais-Smale condition on paths coming from 
monopoles on $[0, \infty) \times T^2$. 
\end{proof} 
 
Since $\chi(T^2)\backslash U_\Theta$ is compact, we can set 
$ 
\delta = \frac 12 min \{ \mu_{a_\infty} | a_\infty \in 
\chi(T^2)\backslash U_\Theta \}. 
$ 
Then, when restricted to the cylindrical 
end,  any Seiberg-Witten monopole on $V$ with finite energy  
and with asymptotic 
limit in $\chi(T^2)\backslash U_\Theta $ 
has an exponential decay at a rate at least $\delta$. 
 
In order to  prove that the boundary value 
map (\ref{boundary-T^2}) is well-defined  and continuous, we need 
to resort to  the ``finite energy implies 
finite length'' principle of L. Simon(\cite{Simon}  
(see Corollary 4.2.5 in \cite{MMR}). 
 
\begin{Rem} 
Given that $f$ is a real analytic function, the work of L. Simon 
as explained in \cite{MMR} 
can be employed to prove a more general {\L}ojasiewicz inequality for $f$ at 
any critical point in $\chi (T^2)$. Let $\gamma (s)$ be a flow line of $f$, 
corresponding to an irreducible 
finite energy solution of the Seiberg-Witten equations on $T^2\times  
[0, \infty)$. Then, 
there exist constants $0< b \le 1$ and $0< c \le \disp{\frac 
12}$ such that, when $s>R >> 1$, we have 
\ba 
\begin{array}{l} 
\inf_{a_\infty \in \chi(T^2)} \|\gamma (s) - a_\infty \|_{L^2} \le 
(\|\nab f(x(s))\|_{L^2})^b, \\[2mm] 
|f(\gamma (s))|^{1-c} \le \|\nab f(\gamma (s))\|_{L^2}. 
\end{array}\label{Loj-f} 
\na 
At the smooth critical points in $\chi (T^2)$, the 
{\L}ojasiewicz inequalities have the best exponents $b = 1$ 
and $c = \frac 12$. 
The direct consequence of these {\L}ojasiewicz inequalities is the following 
finite length result for flow lines: 
\[ 
\int_{s_1}^{s_2} \|\disp{\frac {\d \gamma (s)}{\d s}}\|_{L^2} ds \le 
\disp{\frac 4c} | f(\gamma (s_1))- f(\gamma (s_2))|^c. 
\] 
\label{boundvalue} 
\end{Rem}

Now we have a setting analogous to the key results in
\cite{MMR} (pages 60-70) in our situation. The 
arguments in \cite{MMR}, adapt to the present context, hence imply 
that the boundary value map (\ref{boundary-T^2}) is well-defined and 
continuous as a map 
$$ \partial_\infty:\qquad  \M^*_{V}
\to \chi(T^2), $$ 
In the next subsection, we will study the local properties of the map 
$\partial_\infty$ around the singular point $\Theta$ and the structure
of $\M^*_{V}$. We remind the reader that we have established the exponential
decay property of the monopoles in $\M^*_{V}$
with asymptotic limits away from $\Theta$.

\subsection{Local structure of moduli space of irreducible monopoles} 
 
Let $U_\Theta$ be a small neighbourhood of $\Theta$ in $\chi (T^2)$. 
In this section, we will study the local structure of the  
moduli spaces $\M_V^*\backslash 
(\bar \d_\infty)^{-1}(U_\Theta)$ and $(\bar \d_\infty)^{-1}(U_\Theta)$.  
Here $\bar \d_\infty $ is the 
composition of $\d_\infty$ and $\pi$ in (\ref{boundary-T^2}). 
 
For the structure of the moduli space $\M_V^*\backslash 
(\bar \d_\infty)^{-1}(U_\Theta)$,  
the exponential decay property implies that we can introduce weighted 
Sobolev norms in order to study the Fredholm theory of the 
linearization of the equations.   
With $\delta $ as in the previous subsection, we define the space 
\ba 
{\scr A}_{V, T^2} = \left\{ (A, \psi ) \left| \begin{array}{l} 
(1) \hbox{ $A$ is an extended $L^2_{2, \delta}$-connection on 
$det(W)$}\\ 
(2) \hbox{ $\psi $ is an $L^2_{2, \delta}$-spinor on $W$} 
\end{array}\right. \right\}, 
\label{config-cyl-T2} 
\na 
where extended $L^2_{2, \delta}$-connection means that there exists 
an imaginary-valued  harmonic 1-form $A_\infty$ in $H^1(T^2, i\R)$ such 
that 
$A-A_\infty$ is an $L^2_{2, \delta}$-connection on $det (W)$, where 
$L^2_{2, \delta}$ denotes the Sobolev norm with weight as in 
\cite{LM}. To be precise, 
we choose the weight function $e_\delta(t)= e^{\tilde\delta(t)/2}$, 
where $\tilde\delta(t)$ is a smooth 
function with bounded derivatives, such that 
$\tilde\delta (t)\equiv -\delta t$ for $t\leq -1$ and $\tilde\delta 
(t)\equiv \delta t$ for $t\geq 1$, and for some fixed positive number
$\delta$ defined as 
\[
\delta = \frac 12 min \{ \mu_{a_\infty} | a_\infty \in \chi(T^2) \backslash
U_{\Theta}\}.
\] 
The $L^2_{k,\delta}$ norm is defined as $\| f \|_{2,k,\delta}= \| 
e_\delta f \|_{2,k}$. The weight $e_\delta$ imposes an exponential decay 
as an asymptotic condition along the cylinder.  We define 
the gauge group $\G_{V, T^2}$ to be the $L^2_{3, loc }$--gauge 
transformations such that there exists $g_\infty \in U(1)$  
with $g_\infty^{-1}g-1$ 
an $L^2_{3,\delta}$-gauge transformation.

Assume that $x= (A, \psi) \in {\scr A}_{V, T^2} $ 
is an irreducible ($\psi \ne 0$) perturbed  Seiberg-Witten monopole on
$V$ with finite energy, where the perturbation is in the form of
Section 2 with compact support. Then from the results of previous
subsection, we  
 can assume further that $A_\infty$  represents a flat connection $a_\infty$ 
in $\chi (T^2)\backslash U_\Theta$. Then the irreducible part of the 
fiber $(\bar\d_\infty)^{-1}(a_\infty)$ has a deformation complex 
\ba 
0\to \Lambda^0_{L^2_{3,\delta}}(V, i\R)\stackrel{G}{\to}   
\Lambda^1_{L^2_{2,\delta}}(V, i\R) 
\oplus L^2_{2, \delta}(W) \stackrel{L}{\to} 
\Lambda^1_{L^2_{1,\delta}}(V, i\R)  
\oplus L^2_{1,\delta}(W) 
\label{deform-fiber} 
\na 
where $G$ is the map which gives the infinitesimal gauge 
transformations: 
\[ G\mid_{(A,\psi)}(f)= (-df,f\psi) \] 
and $L$ is the 
linearization  
\ba 
L_{A,\psi} (\alpha,\phi)= \left\{ 
\begin{array}{c} 
      *d\alpha -\sigma(\psi,\phi)\\[2mm]
      \dirac_{A}\phi +\frac 12 \alpha. \psi, 
\end{array}\right. 
\label{L} 
\na 
of the perturbed Seiberg-Witten equations (\ref{SW-YT2}) on $V$. 
We can assemble the deformation complex (\ref{L}) into the 
following operator:
\ba
(G^*_\delta, L):\qquad 
\Lambda^1_{L^2_{2,\delta}}(V, i\R)
\oplus L^2_{2, \delta}(W) \to \Lambda^0_{L^2_{1,\delta}}(V, i\R)
\oplus \Lambda^1_{L^2_{1,\delta}}(V, i\R)
\oplus L^2_{1,\delta}(W)
\label{G-L}
\na
where $(G^*_\delta, L)(\alpha, \phi)$ is given by
\[
(e_{\delta}^{-1} d^* (e_\delta \alpha) + iIm \langle  \psi, \phi\rangle ,
 *d\alpha -\sigma(\psi,\phi), \dirac_{A}\phi +\frac 12 \alpha. \psi).
\]
With the choice of $e_\delta$ as in the previous section, $(G^*_\delta, L)$
is a Fredholm operator of index $0$. 
 
The deformation complex for the moduli space 
$\M_V \backslash (\bar\d_\infty)^{-1} (U_\Theta)$ is 
given by 
\ba 
0\to T_{id}\G_{V, T^2} \stackrel{G}{\to} T_x {\scr A}_{V, T^2} 
\stackrel{L}{\to}  
 \Lambda^1_{L^2_{1,\delta}}(V, i\R) 
\oplus L^2_{1,\delta}(W). 
\label{deform-M} 
\na 
 
These two complexes are related by the fact that (\ref{deform-fiber}) is 
a sub-complex of  (\ref{deform-M}) with the quotient 
complex 
\[ 
0\to Lie (Stab (a_\infty) ) \stackrel{0}{\to} H^1(T^2, i\R) \to 0. 
\] 
Therefore, the virtual dimension of 
$\M_V \backslash (\bar\d_\infty)^{-1}(U_\Theta)$ at $x=(A,\psi)$ is 
$$\dim (\bar\d_\infty)^{-1}(a_\infty) 
+ \dim H^1(T^2, i\R) - \dim Stab (a_\infty)= 1,$$ 
where $\dim (\bar\d_\infty)^{-1}(a_\infty)$ is the virtual dimension 
of the fiber.

\begin{The}  Fix an open set $U$ in $V - (T^2 \times [0, \infty))$. 
There exists a Baire  set ${\cal P}_{0}$ of perturbations 
$\mu$  on $V$ with compact supports in $U$,  
such that  the perturbed Seiberg-Witten moduli space 
$\M_V^*\backslash (\bar\d_\infty)^{-1}(U_\Theta)$ is a smooth, 
oriented manifold of dimension $1$. Moreover, 
\[\d_\infty:
\M_V^*\backslash (\bar\d_\infty)^{-1}(U_\Theta) \to \chi (T^2)
\]
 is an immersion and transversal to any given immersed curves
in $\chi (T^2)$  
\label{smooth-dim1} 
\end{The} 
 
\begin{proof} 
The transversality argument is the same as in the closed 
case, see the proof of Proposition \ref{3dim-compactperturbation}, 
namely, we  look at the deformation complex (\ref{deform-M}) for 
the parametrized moduli space $\M_{V, Z(U, i\R)}$ 
to get the transversality 
for the parametrized moduli space $\M_{V,  Z(U, i\R)}$. We then apply 
the infinite dimensional version of 
Morse-Smale theory to the projection 
$\M_{V,  Z(U, i\R)} \to  Z(U, i\R)$, we obtain that, for $\mu$ in a
Baire space ${\cal P}_{0}\subset  Z(U, i\R)$, the moduli space   
$\M_{V,\mu}^* \backslash (\bar\d_\infty)^{-1}(U_\Theta)$ is a   
smooth  manifold of dimension given by virtual dimension 
calculated as above.  

We first show that for a generic perturbation $\mu$ (a co-closed imaginary
valued 1-form with compact support in $U$), the
map $(G^*_\delta, L)$ as given by (\ref{G-L}) is
surjective. At an irreducible monopole $[A, \psi]$ in 
$(\bar\partial_\infty )^{-1}(a_\infty)$ for $a_\infty \in 
\chi (T^2) \backslash U_\Theta$, we will show that
\[
L: \qquad Ker G^*_\delta \times Z(U, i\R) \rightarrow Ker G^*_\delta
\]
is surjective. Suppose that $(\alpha_1, \phi_1)$ is $L^2_\delta$-orthogonal
to the image of the above map, then 
$(\alpha, \phi) = e_{2\delta}(t) (\alpha_1, \phi_1)$
is
$L^2$-orthogonal to the image of the above map, hence, $(\alpha, \phi)$
is in $L^2_{1, -\delta}$ and 
satisfies the equations (\ref{3d-surj}) as in the proof of
Proposition \ref{3dim-compactperturbation}. 
Hence, there is a real valued function $f$ on $V$ (with infinite
cylindrical end) such that
$\phi = if \psi$, $\alpha = -2i df$
and
\[
2d^*d f + |\psi|^2 f =0.\]

$df \in L^2_{1, -\delta}$ implies that on $T^2\times [0, \infty)$,
$\displaystyle {\frac {\partial f}{\partial t}}$ is in $L^2_{-\delta}$, 
then by Cauchy-Schwartz inequality
\[
|f (t) -f (0)|^2 \le \displaystyle{
\int_0^t e^{2\delta s}ds \int_0^t |e^{-\delta s}\frac {\partial f}{\partial t}
(s)|^2 ds }\]
this implies that for $T>>0$,
\ba
\displaystyle{
\int_{\partial V(T) } |f| ^2 }\le C_0 e^{2T\delta} \|
\displaystyle {\frac {\partial f}{\partial t}}\|^2_{L^2_{-\delta}}
\le Ce^{2T\delta},
\label{inequality-1}
\na
for some constants  $C_0, C$.
$[A, \psi] \in (\bar\partial_\infty )^{-1}(a_\infty)$, whose
asymptotic behaviour has been studied in the previous subsection, we see
that there exist gauge representatives $(A, \psi)$ and $a_\infty$
of $[A, \psi]$ and $[a_\infty]$, so that
 $(A, \psi)$ decays to $(a_\infty, 0)$ exponentially at the rate
at least $\displaystyle{ \frac {\mu_{a_{\infty}}}{2}}$, where 
$\mu_{a_{\infty}}$ is
the smallest absolute value of the non-zero eigenvalue
of $Q_{a_{\infty}}$ (Cf. (\ref{Hessian-f})). 
On $T^2\times [0, \infty)$, write 
$\alpha = \pi^*(a_{\infty})+  
\alpha^1 + \alpha^0 dt$ with $\alpha^i \in \Omega^i (T^2, i\R)$,  
using the analysis in Appendix of \cite{Mrowka}, we get
\[
\|\alpha^0 \|_{C^0(T^2\times [T, T+1])} \le 
C_1 e^{-T\displaystyle{ \frac {\mu_{a_{\infty}}}{2}}},\] from
some constant $C_1$.
As $-2 i\displaystyle{\frac {\partial f}{\partial t} } = \alpha^0 $, 
we obtain   
\ba
\left|\displaystyle{\frac {\partial f}{\partial t} }|_{\partial
V(T)}\right| \le C_2 e^{-T\displaystyle{ \frac {\mu_{a_{\infty}}}{2}}},
\label{inequality-2}
\na
for some constant $C_2$. From two inequalities (\ref{inequality-1}) and
(\ref{inequality-2}), and  note that $\delta \le \mu_{a_{\infty}}/2$, we have
\[
\begin{array}{lll}
&& \left|\displaystyle{\int_{V}} d^*d (f^2) dvol\right| \\[2mm]
&\le & \lim_{T\to\infty} \displaystyle{\int_{\partial V(T)} 
2 |f| \cdot |\frac{\partial f}{\partial t} }| dvol \\[2mm]
&\le& \lim_{T\to\infty} 2 \displaystyle{
(\int_{\partial V(T)} |f|^2 dvol )^{\frac 12} (\int_{\partial V(T)} |
\frac{\partial f}{\partial t} |^2 dvol )^{\frac 12}} \\[2mm]
&\le& \lim_{T\to\infty} 2 CC_1 e^{ T 
(\delta-\displaystyle{\mu_{a_{\infty}}/2} )}=0.
\end{array}
\]
Now multiply $2d^*df + |\psi|^2 f =0$ by $f$ and integrate it by parts,
we get $f=0$, hence $(\alpha_1, \phi_1) =0$.

The proof of that$\partial_\infty$ is an immersion
and transversal to any given immersed curves follows from the 
Sard-Smale theorem.

An orientation of $\M_V$ is obtained from a trivialization of the 
determinant line bundle of the assembled operator of the deformation 
complex (\ref{deform-M}). The trivialization of 
the determinant line bundle of the 
complex of  (\ref{deform-fiber}) is obtained from the 
orientation of $H^0_\delta(V)\oplus H^1_\delta(V)$, 
the cohomology groups of $\delta$-decaying forms. 
In fact, we can deform the operator 
$H_{(A,\psi)}$ with a homotopy $\epsilon \psi$, $\epsilon\in [0,1]$. 
The asymptotic operator $Q_{[a_\infty,0,0]}$ is preserved 
in the deformation. Thus, if the weight $\delta$ is chosen in such a 
way that $\delta/2$ is not in the spectrum of $Q_{[a_\infty,0,0]}$, 
then (\cite{LM}, \cite{MMR} Lemma 8.3.1) we can ensure that the 
operator $H_{(A,\epsilon \psi)}$ is Fredholm, for all 
$\epsilon\in [0,1]$. Since the Dirac operator is complex linear and it 
preserves the orientation induced by the complex structure on the 
spinor bundle, a trivialization of the determinant line bundle at 
$\epsilon= 0$ is obtained by the orientation of 
$H^0_\delta(V)\oplus H^1_\delta(V)$. This in turn determines a 
trivialization of the determinant line for $\epsilon= 1$, hence an 
orientation of $\M_V$. 
\end{proof} 
 
Similar results were obtained by \cite{Chen} \cite{Lim}.

Now we need to understand the local structure of $\M^*_V$ around 
$(\bar\d_\infty)^{-1}(\Theta)$.  
The center manifold technique 
developed in \cite{MMR} is a useful model to study the structure of 
$(\bar\d_\infty)^{-1}(U_\Theta)$.  
 
We briefly recall a few facts about center manifolds \cite{MMR}. In general, 
suppose we are given a system of the form 
\ba 
\dot{x}=Q x + N(x),  
\label{aut-diff} 
\na 
with $Q$ a linear operator acting on a Hilbert space ${\cal 
X}$. Assume we also  
have the decomposition ${\cal X}={\cal X}_h^+ \oplus {\cal 
X}_c \oplus {\cal X}_h^-$ determined by the positive, negative, and zero  
spectrum of the operator $Q$. Let ${\cal X}_h={\cal X}_h^+ \oplus 
{\cal X}_h^-$, and consider the projections $\pi_c : {\cal X}\to {\cal 
X}_c$ and $\pi_h: {\cal X}\to {\cal X}_h$. We denote by $Q^\pm_h$ and 
$Q_c$ the induced operators on ${\cal X}_h^\pm$ and ${\cal X}_c$. By 
construction $Q_c$ is trivial.  
The evolution semigroups $e^{-s Q_h^+}$ and $e^{s Q_h^-}$, for 
$s\geq R_0>0$, satisfy  
\ba 
\sup_{s\geq R_0}\max \{ e^{\delta s}\| e^{-s Q_h^+} \| , 
e^{-\delta s}\|  e^{s Q_h^-} \| \} \leq C,  
\label{bound-evol} 
\na 
for some constant $C >0$. This follows from the bound 
$$ \frac{1}{2} \inf \{ | \lambda | \ | \ \lambda\in spec(Q_\Theta), 
\lambda\neq 0 \} = \delta > 0. $$ 
 
The {\bf center manifold theorem} (in \cite{MMR}) 
states that there  
exists a  map $\varphi: {\cal X}_c \to {\cal X}_h$ that vanishes  
to second order at the origin, and such that an 
element $\tilde x(s)$ is a solution of (\ref{aut-diff}) if and only if  
the projection $\pi_c \tilde x(s)$ is a solution of the equation 
\ba 
\dot{x_c}=\pi_c N(x_c + \varphi(x_c)). 
\label{aut-diff-c} 
\na 
The {\em center manifold} ${\cal H}$ is defined as 
$ {\cal H}=\{ x_c + \varphi(x_c) |  x_c \in {\cal X}_c \}$.

We now describe explicitly the center manifold and the stable set for the 
unperturbed equations (\ref{SW-T2}). In this case, 
we are considering the operator $Q_\Theta$, the Hessian of 
the functional $f$ at the degenerate critical point $\Theta$.  
The center manifold ${\cal H}_\Theta$ for the functional 
$f$ at the degenerate critical point $\Theta$ is a 
$C^2$-manifold which is invariant under the gradient flow of $f$, 
contains a small neighbourhood $U_\Theta$ of $\Theta$, and has tangent 
space at $\Theta$ given by 
\[ 
{\cal H}^1_\Theta = H^1(T^2, i\R) \oplus ker \bar\d_\Theta 
\oplus ker \bar\d_\Theta^*  \cong \CC^3. 
\]

\begin{Lem} At every point $x=(a, \alpha, \beta )\in {\cal H}^1_\Theta$, 
the gradient vector $\nab f (a, \alpha, \beta)$ is tangent to 
${\cal H}^1_\Theta$, hence ${\cal H}^1_\Theta$ is a center manifold 
of $f$ around $\Theta$. 
\end{Lem} 
 
\begin{proof} Using the natural complex structure on $T^2$, we can 
identify 
${\cal H}^1_\Theta$ as the space of constant sections of 
\[ 
\Lambda^1(T^2, i\R) \oplus \Lambda^{0,0}(T^2, \CC) \oplus 
\Lambda^{0,1}(T^2, \CC).\] 
For $(a, \alpha, \beta )\in {\cal H}^1_\Theta$, we have 
$ 
\nab f (a, \alpha, \beta ) = (-i\bar \alpha \beta, -i\bar\d^*_a\beta, 
i\bar\d_a \alpha), 
$ 
which is a constant section. 
Take $(z_1, z_2, z_3)$ as the coordinates on ${\cal H}^1_\Theta 
\cong \CC^3$, we have 
\[ 
\nab f(z_1, z_2, z_3) = \bigl(-\bar z_2 z_3, -\frac{\bar z_1z_3}{2}, 
-\frac{z_1z_2}{2}\Bigr). 
\] 
\end{proof} 
 
The downward gradient flow of $f$ on ${\cal H}^1_\Theta$ is given by 
\ba\left\{ 
\begin{array}{l} 
\disp{\frac {\d z_1}{\d s}} = \bar z_2 z_3, \\[2mm] 
\disp{\frac {\d z_2}{\d s}} = \frac{\bar z_1 z_3}{2}, \\[2mm] 
\disp{\frac {\d z_3}{\d s}} =  \frac{z_1z_2}{2}. 
\end{array} 
\right. 
\label{gradient-CM} 
\na 
Note that this gradient flow is invariant under the $U(1)$-action 
(the constant gauge transformation): 
\[ 
(z_1, z_2, z_3) \stackrel{e^{i\alpha}\in U(1) }{\to} (z_1, e^{i\alpha}z_2, 
e^{i\alpha}z_3).\] 
 
\begin{Lem}  The quantities $|z_2|^2 -|z_3|^2$, 
$|z_1|^2 - |z_2|^2 -|z_3|^2$ and $Im(z_1z_2\bar z_3)$ are preserved 
under the 
gradient flow on ${\cal H}_\Theta$. 
\label{preserve-flow} 
\end{Lem} 
\begin{proof} This is a direct calculation using the gradient flow 
equations (\ref{gradient-CM}). \end{proof} 
 
The stable set of $(a_\infty, 0)\in {\cal H}_\Theta$ is defined to be 
\[ 
{\cal S}_{a_\infty} = \{ x\in {\cal H}_\Theta\  \hbox{such that the 
flowline 
of (\ref{gradient-CM}) starting at $x$ converges to $(a_\infty, 0)$} 
\}. 
\] 
The stable set ${\cal S}$ of $f$ in ${\cal H}_\Theta$ is the union of these 
${\cal S}_{a_\infty}$, for $(a_\infty, 0)\in {\cal H}_\Theta$. 
 
\begin{Lem} Let $(a, \phi)= (z_1,  z_2, z_3)\in {\cal H}_\Theta$. 
Then $ (a, \phi) \in {\cal S}$ if and only if we have 
\ba 
\left\{\begin{array}{l} 
|z_2|^2 -|z_3|^2 = 0, \\[2mm] 
|z_1|^2 - |z_2|^2 -|z_3|^2 = |a_\infty|^2 \ge 0, \\[2mm] 
Im(z_1z_2\bar z_3) = 0. 
\end{array} 
\right. 
\label{stable-general} 
\na 
In particular, $(a, \phi)= (z_1,  z_2, z_3)\in {\cal S}_\Theta$ 
(the stable set of the point $\Theta$) if and only if 
\ba 
\left\{\begin{array}{l} 
|z_2|^2 -|z_3|^2 = 0, \\[2mm] 
|z_1|^2 - |z_2|^2 -|z_3|^2 = 0, \\[2mm] 
Im(z_1 z_2\bar z_3) = 0. 
\end{array} 
\right. 
\label{stable-Theta} 
\na 
These equations describe a cone  over a torus $T^2$ with vertex at $\Theta$. 
Furthermore,  ${\cal S}\backslash \{(\Theta, 0)\}$ is a 
4-dimensional manifold with boundary 
 ${\cal S}_\Theta\backslash \{(\Theta, 0)\}$. 
\end{Lem} 
\begin{proof} It follows from Lemma \ref{preserve-flow} that $(a, \phi) 
\in {\cal S}$ converges to some $(a_\infty, 0)$ as $t \to \infty$. 
The equations (\ref{stable-Theta}) 
define a torus over $\Theta$ and, as $|a_\infty|^2 \to 0$, 
points defined by (\ref{stable-general}) approach points in ${\cal 
S}_\Theta$. 
\end{proof}

As in \cite{MMR}, the restriction of the gradient  
flow of $f$ to the center manifold 
provides a model for the structure of the space of flows with 
asymptotic values in a small neighbourhood $U_\Theta$ of  
 $\Theta$, in the following sense. 
Given a point $x$ in ${\cal H}_\Theta$, the {\em stable set} ${\cal 
S}_x$ at $x$ is defined as 
$$ {\cal S}_x= \{ y \in {\cal H}_\Theta\  \hbox{such that the flowline 
starting at $y$ converges to $x$} 
\}.$$ 
The stable set ${\cal S}= \cup_{x\in{\cal H}_\Theta}{\cal S}_x$ defines 
a 
refinement of the boundary value map, as described in the following 
commutative diagram: 
\ba 
\diagram 
(\bar\d_\infty)^{-1} (U_\Theta) \rrto^{\bar\d_\infty} \drto_\Upsilon 
& & U_\Theta \\ 
& {\cal S} \urto^{\Gamma} & \\ 
\enddiagram 
\label{refine-map} 
\na 
Here $\Upsilon$  is a map defined by taking a flow line on the 
stable set ${\cal S}$ that is exponentially close to a monopole in 
$(\bar\d_\infty)^{-1} (U_\Theta)$. The map $\Gamma$ is the limit 
value map under the flow line of $f$ on ${\cal S}$. 
 
The  results of \cite{MMR} show that if the projection $\pi_c \tilde x(s)$ of a flow line $\tilde 
x(s)$ satisfies an estimate  
\ba 
\| \partial_s \pi_c \tilde x(s) - \nabla f (\pi_c \tilde x(s)) \| < 
C e^{-\delta s}, 
\label{estMMR} 
\na 
for all $s\geq R_0$, then there exists a unique flow line $x_c(s)$ in 
the center manifold ${\cal H}_\Theta$ that is exponentially close to $\pi_c 
\tilde x(s)$ for large $s \geq R_0$, with the same exponent $\delta$ 
determined by the smallest absolute value of the non-zero 
eigenvalues of $Q_\Theta$. Moreover, for a flow line $\tilde 
x(s)$ satisfying 
$$ \| \pi_c \tilde x(s) - a_\infty \|_{L^2_2(T^2\times \{ s \})} + \| 
\pi_h \tilde x(s) \|_{L^2_2(T^2\times \{ s \})} \leq C $$ 
for all $s\geq R_0$, the projection $\pi_h \tilde x(s)$ is 
exponentially small for large $s$, with exponent $\delta$. 
The condition (\ref{estMMR})  follows from our explicit 
construction of the center manifold. This shows that every flow line 
in $(\bar\d_\infty)^{-1} (U_\Theta)$ is exponentially close to a 
flow line in the center manifold. Thus, the refinement $\Upsilon$ of 
the boundary map is well defined and continuous. 
 
The results of the previous discussion and the arguments in \cite{MMR} 
(Page 82-100)  imply  
the following structure theorem for our moduli space 
$\M_{V}$ near $\bigl(\bar\d_\infty\bigr)^{-1} (\Theta )$.    
 
\begin{The}  
\label{3-d-cyl}  
Fix a metric $g$ and perturbation $P\in {\cal P}$ 
 as in Theorem \ref{smooth-dim1}. 
Let ${\cal K} \subset \M^*_{V, P}$ be defined as ${\cal K}= 
\bigl(\bar\d_\infty\bigr)^{-1} (\Theta )$, and let ${\cal K}'$ denote 
the subset $\Upsilon^{-1}(\Theta ,0, 0)$, with $\Upsilon$ defined as 
in (\ref{refine-map}). Then, generically, the following holds. 
 
(1) ${\cal K}'$ is empty and ${\cal K}$ consists of only finitely 
many points. 
 
(2) There is a neighbourhood $U_\Theta$ of 
$\Theta$ in $\chi(T^2)$, such that the following holds. The moduli 
space  
$ \M^*_{V}\cap \bigl(\bar\d_\infty\bigr)^{-1}(U_\Theta )$  
is a smooth manifold of dimension $1$, with boundary ${\cal K} = 
\bigl(\bar\d_\infty\bigr)^{-1} (\Theta )$.  
\end{The} 
 
\begin{proof} 
By the center manifold theorem in \cite{MMR}, the restriction of any 
finite energy monopole  
$$[A, \psi] \in \M_{V, \mu}(U_\Theta) = 
\bigl(\bar\d_\infty\bigr)^{-1} (U_\Theta)$$ 
to the tube $T^2 \times [T_0, \infty)$ (for a fixed large $T_0$)  
is exponentially close to a flow line in the center manifold starting 
from the point $\Upsilon ([A, \psi])$ given by the  
refinement boundary map 
(\ref{refine-map}). The exponential weight is at least a half of the 
smallest absolute value of the non-zero eigenvalues of 
$(\bar\d_\Theta + \bar\d_\Theta^*)$.  Theorem \ref{smooth-dim1} 
shows that, for a generic choice of the 
perturbation, the moduli space $\M_{V}^*(U_\Theta)$ 
is a smooth manifold of dimension $1$,  away from $\Upsilon^{-1}({\cal 
S}_\Theta)$. From the analysis of the center manifold theorem, since  
${\cal S}\backslash \{(\Theta, 0)\}$ is a 4-manifold with  
boundary ${\cal S}_\Theta \backslash \{(\Theta, 0)\} $, 
we know that generically ${\cal K}'$, if non-empty, is a smooth manifold  
of dimension given by the virtual dimension: $ dim \M_{V}^* 
-4 =-3$, so ${\cal K}'$ must be empty and $\M_{V}^*(U_\Theta)$  
is a smooth oriented 1-dimensional manifold with boundary ${\cal K}= 
\bigl(\bar\d_\infty\bigr)^{-1} (\Theta )$.  
\end{proof} 
 
One useful observation that we can derive directly from the analysis of 
the center manifold is the following estimate of the rate of decay of 
solutions approaching the singular point $\Theta$. 

\begin{Rem} 
Let $x(s)=(a(s),\alpha(s),\beta(s))$ be an irreducible finite energy 
solution of the Seiberg-Witten equations on $V$, with asymptotic 
value $\Theta$, that is, $[x]\in \partial_\infty^{-1}(\Theta)$. Then 
the rate of decay in the $s\to\infty$ direction is polynomial with 
$$ \| (a(s)-\Theta,\alpha(s),\beta(s)) \|_{L^2(T^2\times \{ s\})} 
\sim \frac{1}{s}. $$ 
\label{1/s-decay} 
\end{Rem}

\subsection{Proof of Theorem \ref{structure}} 
 
  From the discussions in the previous subsection, in order to
complete our analysis of the structure of the moduli space $\M_V$, 
we only need to prove the following result.

\begin{Lem} 
$\M_V^*$ is compact except for finitely many open ends limiting to 
$\chi (V)$, the reducible moduli space of $V$ and, generically, 
$ \bar\d_\infty (\M_V^*)$ can be made transversal   
at any interior points to any given finite 
set of curves in 
$\chi(T^2, V)$. 
\label{lastMV}
\end{Lem}
   
\proof
We first analyse the set of reducible solutions of the monopole equations on 
$V$. The reducible moduli 
space $\M_V^{red}$ can be identified with the space 
$\chi (V)$ of deformed flat connections over $V$, modulo gauge 
transformations, which is diffeomorphic to a circle.  
The asymptotic value map $\partial_\infty$  
is simply the restriction 
map, which is an embedding 
$$ \partial_\infty : \M_V^{red} =  \chi(V) \hookrightarrow \chi(T^2, V). $$ 
 
Let $\chi(V)\hookrightarrow 
\chi(T^2,V)$ be the circle of reducibles on $V$ modulo 
gauge equivalence, embedded via the restriction map inside the 
cylinder $\chi(T^2,V)$. 
Fix a smooth parameterization $a(t)$ of $\chi(V)$, 
consider the family of 
Dirac operators $\dirac_{a(t)}$ on $V$, twisted with the connection 
$a(t)$. We can perturb $\chi(V)$  
such that $\chi(V)\hookrightarrow 
\chi(T^2,V)$ is  away 
from a small neighbourhood of the singular point $\Theta$, then 
we know that the Dirac operator $\bar\d_{a(t)}+\bar\d^*_{a(t)}$ on $T^2$ has 
trivial kernel. The 3-dimensional Dirac operator $\dirac_{a(t)}$ 
($a(t) \in \chi (V)$ ) on $V$ may acquire a non-trivial kernel,
however, this only happens at finitely many points on $\chi (V)$, for  
a generic perturbation in ${\cal P}_{0}$ (cf.~ \cite{MW} \S 7).  
We show that, if these occur, then  
the irreducible set $\M_V^*$ has an open end limiting to  
such points. If the irreducible set $\M_V^*$  
has an open end limiting to the reducible set $\chi (V)$, then  
the 3-dimensional Dirac operator $\dirac_{a(t)}$ has a non-trivial
kernel: this can be seen by studying the linearization of the 
spinor part of the Seiberg-Witten equations. On the other hand,  
suppose that there is a point  
$a(t_0)$ on $\chi (V)$, where the  
operator $\dirac_{a(t_0)}$ acquires a non-trivial kernel.  
 We can proceed as in Section 7.3 of \cite{MW} 
to analyse the local model of the moduli space $\M_V = \M_V^* \cup \chi (V)$  
in a neighbourhood of $[a(t_0),0]$, which shows that there  
exists an open end limiting to $a(t_0)$.  
 
Thus, the rest of the proof of compactness for $\M_V = \M_V^* \cup
\chi (V)$ is now reduced to the (by now standard) proof of compactness
for Seiberg-Witten moduli spaces \cite{KM1}
\cite{Lim2} \cite{Mor}.  
Transversality at the interior  
of the $\bar\d_\infty (\M^*_V)$ to any given finite set of  
curves in $\chi(T^2, V)$ can also be achieved by a 
generic choice of perturbation on $V$. 

\endproof

Thus, we have completed the proof of the structure theorem for $\M_V$
(Theorem \ref{structure}). 

\section{Gluing  of 3-dimensional monopoles}

Now we begin to discuss the gluing theory. 
Suppose that $V(r) = V - (T^2 \times [r, \infty))$ lies in  
a closed 3-manifold $Z$ such that $T^2$ splits $Z$ into two components 
(for example, $V(r) \cup _{T^2} \nu(K)$). We identify 
the solutions of the Seiberg-Witten equations on $V$ differing only 
by those gauge transformations $V$ which can be 
extended to $Z$, and denote the resulting 
moduli space by $\M^*_{V,Z}$. Then the boundary value map 
in Theorem \ref{structure} has a refinement: 
\[ 
\M^*_{V,Z} \longrightarrow \chi(T^2, Z) 
\] 
where  the notation 
$\chi (T^2, Z)$ indicates the moduli space of flat connections on a trivial 
line bundle over $T^2$ modulo the gauge transformations on $T^2$ which can be 
extended to $Z$. 
This gives a refined boundary value map and 
the moduli spaces $\M^*_{V,Z}$  enjoy all the properties described in  
Theorem \ref{structure} for $\M^*_V$. 

Assume that $Z = V(r) \cup_{T^2} \nu(K)$ where $\nu(K)$ is 
a tubular neighbourhood of a knot $K$ in $Z$. 
We denote by $\chi (\nu(K) , Z)$ 
the moduli space of flat connections 
on $\nu(K)$ modulo the gauge transformations on $\nu(K)$ which can be 
extended to $Z$. There is a natural map 
$\chi (\nu(K) , Z) \hookrightarrow \chi (T^2, Z)$,
which realizes $\chi(\nu(K) , Z)$ as a line in the affine space
$\chi (T^2, Z)$. 
 
Thus we can define the following fiber product: 
\ba 
\M^*_{V,Z}\times_{\chi(T^2, Z)} \chi (\nu(K), Z). 
\label{fiber-product} 
\na 
This is the main object in the gluing Theorem \ref{Gluing}.
We shall present the argument for the case of the homology sphere
$Y$. The argument is analogous in the case of $Y_1$ and, up to minor
modifications that we shall point out, in the case of $Y_0$ as well.

Consider a tubular
neighborhood $\nu(K)\subset Z$ endowed with a metric with sufficiently
large positive curvature inside $\nu(K)$ and flat near the boundary.
When stretching the neck in $Z(r)$, using the standard
pointwise estimate on the spinor for Seiberg-Witten monopoles we can
ensure that, on $\nu(K)$ endowed with an infinite cylindrical end, 
the only finite energy solutions of the unperturbed Seiberg-Witten
equations are reducibles (with vanishing spinor part). 
Modulo gauge transformations,   
these correspond to the moduli space of flat connections 
on $\nu (K)$.  
In Lemma \ref{Liviu-metrics}, we will show that, if we choose such a
metric for $\nu(K)\subset Y$, it is still possible to have a metric
with the same properties for $\nu(K)\subset Y_1$ and $\nu(K)\subset
Y_0$. 
 
Recall that we have a splitting of 
$Y$ along the torus $T^2$ as $Y= V \cup_{T^2} \nu(K)$, 
with $\partial V= \partial \nu(K)= T^2$. Assume that the metric $g$ on 
$Y$ is the product metric on a small neighbourhood of 
$T^2$, and can be extended to a metric on $\nu(K)$ 
with positive scalar curvature. 
On both $V$ and $\nu(K)$ we consider as underlying $Spin$ 
structure the one induced from the restriction of the trivial $Spin$ 
structure on $Y$. This induces a non-trivial $Spin$ structure on 
$T^2$. The corresponding $\spinc$ structures $\s'$, $\s''$  
on $V$ and $\nu(K)$ have trivial determinant. 
In gluing the $\spinc$ structures $\s'$ and $\s''$  on $V$ and $\nu(K)$ 
we can only obtain the unique trivial $\spinc$ structure 
on  $Y$ since $Y$ is a homology sphere. 
The same holds for $Y_1$. In the case of $Y_0$, 
the gluing of the trivial structures $\s'$ and $\s''$  on $V$ and $\nu(K)$ 
by  gauge transformations along the common boundary $T^2$ 
provides different $\spinc$ structures on $Y_0$, which 
are classified by  
\[ 
H^1(T^2, \Z) /\{Im (H^1(V, \Z), H^1(T^2, \Z)) +  
Im (H^1(\nu(K), \Z), H^1(T^2, \Z))\} \cong H^2(Y_0, \Z). 
\] 
Thus, there are a $\Z$-family 
of $\spinc$ structures corresponding to $H^2(Y_0, \Z)\cong \Z$. 
 
Let $Y(r)= V \cup_{T^2} ([-r,r]\times T^2) \cup_{T^2} \nu(K)$. We can also 
consider the manifolds $V$ and $\nu(K)$ 
with infinite cylindrical ends as 
\[   V \cup_{T^2} ([0, \infty ) \times T^2),\quad \quad 
( (-\infty, 0]  \times T^2) \cup_{T^2} \nu(K). 
\] 
We continue to use the same notation $V$ and $\nu(K)$ for the 
manifolds with infinite cylindrical ends, as we did in the previous
sections. 
 
The proof of the gluing Theorem \ref{Gluing} consists of several steps. 
First, we show that, upon stretching the 
neck $[-r,r]\times T^2$ to infinity, the Seiberg-Witten monopoles 
 $ (A_r,\psi_r)$ on $Y(r)$ approach a 
pair of finite energy solutions $(A',\psi')$, and $(A'',0)$  
on the two manifolds $V$ and $\nu(K)$ with infinite cylindrical
ends. Then we construct a gluing map, under the hypothesis that the 
gluing takes place away from  $\Theta$ in the 
character variety $\chi(T^2)$. 
At the end of this section, we justify the 
assumption that gluing at $\Theta$ can be avoided.

\subsection{Convergence of monopoles on a 3-manifold with a long neck} 
 
We need to introduce some {\em ad hoc} assumptions on the class ${\cal  
P}$ of perturbations for the Seiberg-Witten monopoles on 
$Y(r)$, so that it behaves nicely under the splitting 
$r\to\infty$. We consider perturbations of the monopole equations as 
in (\ref{3dSW}), induced by the perturbations of the
Chern-Simons-Dirac functional. Notice that, if we choose   
a perturbation with compact support on the manifold $V$ with infinite
cylindrical end, this perturbation induces a perturbation on $Y(r)$,
for sufficiently large $r > > r_0$, which is supported inside the knot 
complement in $Y(r)$ (which we still denote by $V$). 

The convergence result we prove in this section depends on 
a uniform pointwise bound on the  
solutions $(A_r,\psi_r)$ in ${\cal M}_{Y(r)}$ which is independent of 
$r$. The argument for the manifold $Y_1(r)$ is the same. The case 
of the manifold $Y_0(r)$ is also analogous, whenever $Y_0$ is endowed 
with a $\spinc$ structure that restricts to the trivial 
$\spinc$ structures on $V$ and 
$\nu(K)$.  

In order to derive the estimates we need, we consider first, for $Y$
a 3--manifold (either without boundary, or with boundary $T^2$) a
functional on the configuration space of $U(1)$--connections and
spinors of the form 
{\small \begin{equation}
{\cal E}_{Y,\mu} (A,\psi)=\int_Y \left( |\nabla_A \psi|^2 + \frac{\kappa}{4}
|\psi|^2  + \frac{1}{2} | F_A |^2 + \frac{1}{2} |
\sigma(\psi,\psi) +\mu |^2\right) \, dv.
\label{Henergy} \end{equation} }
Here we consider compactly supported perturbations $\mu$ of
the form described in Section 2.

\begin{Lem}
If $(A,\psi)$ is a solution of the perturbed SW equations 
$$ (*F_A - \sigma (\psi, \psi) - \mu, \dirac_A\psi ) =(0, 0),
$$ on a compact 3--manifold $Y$ without boundary,
then we obtain
\begin{equation} \label{H1} {\cal E}_{Y,\mu} (A,\psi)= 
 \int_Y F_A \wedge \mu. \end{equation}
If we consider an open submanifold $Z\subset Y$ with boundary
$\partial Z =T^2$, such that the perturbation $\mu$ is
supported away from $\partial Z$, 
then for the functional ${\cal E}_{Z,\mu}$ we have 
\begin{equation} {\cal E}_{Z,\p} (A,\psi)= \int_Z F_A 
\wedge \mu - \int_{\partial Z} \langle \alpha, i \bar\partial_a^* \beta
\rangle, \label{H2} \end{equation}
where we write the connection and spinor as $(a,\alpha,\beta)$ on $T^2
=\partial Z$.
In particular for a cylinder region $Z=T^2 \times [s_0, s_1]$, and 
perturbation term $\mu$ supported away from $Z$, we have
\begin{equation} \label{H3} f(a(s_1),\alpha(s_1),\beta(s_1)) -
f(a(s_0),\alpha(s_0),\beta(s_0)) 
= {\cal E}_{Z,\mu} (A,\psi), \end{equation}
where we write $(A,\psi)$ in the form $(a(s),\alpha(s),\beta(s))$ on
the cylinder $Z$. 
\label{Hmonopoles}
\end{Lem}

\proof
First notice that we have
$$ \int_Y |\dirac_A \psi |^2 \, dv = 
 \int_Y |\nabla_A \psi|^2 + \frac{\kappa}{4}
|\psi|^2 - \frac{1}{2} \langle *F_A \cdot \psi, \psi \rangle . $$
Here the term $\frac{1}{2} \langle * F_A \cdot \psi, \psi \rangle$
can be written as $-F_A \wedge \sigma(\psi,\psi)$. 
We also have
$$ \int_Y | *F_A -\sigma(\psi,\psi)-\mu |^2 \, dv = $$
$$ \int_Y | F_A |^2 + | \sigma(\psi,\psi) + \mu |^2 \, dv 
+2 \int_Y F_A \wedge (\sigma(\psi,\psi) + \mu). $$

Thus, we can rewrite the functional \eqref{Henergy} in the form
$$ {\cal E}_{\mu} (A,\psi)= \int_Y |\dirac_A \psi |^2 +
\frac{1}{2} | *F_A -\sigma(\psi,\psi)-\mu |^2 \, dv
 + \int_Y F_A \wedge \mu. $$ 
The identity \eqref{H1} for a compact manifold then follows. In the
case of  \eqref{H2} for $Z$ 
with $\partial Z=T^2$, see the proof of Lemma \ref{variation1},
the boundary term is the difference of
$$ \int_Z \left( |\dirac_A \psi |^2 - \langle \dirac_A^* \dirac_A \psi, \psi
\rangle \right)\, dv $$
and
$$ \int_Z \left(|\nabla_A \psi|^2 - \langle \nabla_A^* \nabla_A \psi, \psi
\rangle \right)\, dv. $$ 
The last case \eqref{H3} for a cylinder follows, since by the
assumption on 
the perturbation term $\int F_A\wedge \mu$
is trivial, and the boundary terms give the variation of the
functional $f$ along the cylinder. 
\endproof

Notice that the above result allows us to obtain estimates for
the $L^2$ norms of $\psi$, $\nabla_A \psi$, and $F_A$.

\begin{Lem} 
Suppose we are given solutions $(A_r,\psi_r)$ of the perturbed 
Seiberg-Witten equations (\ref{3dSW}) 
on the compact 3-manifold $Y(r)$, with a perturbation $\mu$ supported
in the knot complement $V\subset Y(r)$ for all $r\geq r_0$. Then we have 
pointwise bounds  
\[ | \psi_r (y) |\leq \kappa(Y),  \ \ \ \ \ 
 | F_{A_r} (y) | \leq C (\kappa(Y))^2, \]  
for $y\in Y(r)$, 
where $C, \kappa(Y)$ are constants independent on $r$. 
\label{pointboundLem} 
\end{Lem}

\proof
Consider $\kappa(Y(r))= \max_{y\in Y(r)} \{ -\kappa(y)+C, 0 \}$, where
$\kappa(y)$ is the  
scalar curvature and $C$ is a constant depending only on the 
perturbation $\mu$. Notice that, 
by our assumptions on the choice of the perturbation, we can 
assume that $C$ is independent of $r$.  
The minimum of the scalar curvature also remains constant 
upon stretching the cylinder $T^2\times [-r,r]$, so that we have 
$\kappa(Y(r))= \kappa(Y)$ for all $r>0$. 
 
The Weitzenb\"ock formula provides a uniform bound on the spinors in 
terms of the scalar curvature,
namely at a point $y$ where  
$|\psi_r (y)|$ achieves a maximum we have either $\psi_r(y)= 0$ or 
$ | \psi_r (y) |^2 \leq -\kappa(y)+C. $ 
The pointwise bound for the curvature form $F_{A_r}$ follows from the 
bound on $| \psi_r |$ and from the equations. 
\endproof

Using these pointwise estimates and the results of Lemma
\ref{Hmonopoles}, we obtain $L^2$ and $L^2_1$ estimates. 

\begin{Lem} 
Suppose we are given solutions $(A_r,\psi_r)$ of the perturbed 
Seiberg-Witten equations (\ref{3dSW}) 
on the compact 3-manifold $Y(r)$, with a perturbation $\mu$ supported
in the knot complement $V\subset Y(r)$, for $r\geq r_0$.  

(i) Consider an open submanifold $Z\subset Y(r)$, with $\partial
Z=T^2$ a slice in the product region of $Y(r)$. Then the values
$f(a_r,\alpha_r,\beta_r)$ on $\partial Z$ are 
uniformly bounded in $r\geq r_0$. Here the $(a_r,\alpha_r,\beta_r)$ are
restrictions to $\partial Z$ of the solutions $(A_r,\psi_r)$.

(ii) The total variation of the functional $f$ along a cylinder $Z_r=
T^2\times [-r,r]\subset Y(r)$ is uniformly bounded in $r\geq r_0$.

\label{Hbound}
\end{Lem}

\proof
Applying \eqref{H1} of the previous Lemma together with the
assumptions on the perturbation, we obtain
$$ c \leq - \frac{\kappa(Y(r_0))}{4} \| \psi_r \|^2_{L^2(Y(r_0))} \leq
{\cal E}_{Y(r),\mu}(A_r,\psi_r) =  \int_{Y(r_0)} F_{A_r} \wedge
\mu \leq C', $$ 
with $\kappa(Y(r_0))=\max(-\kappa(y)+C,0)$, for $y\in
Y(r_0)$. We are using the fact that the scalar curvature satisfies
$\kappa\equiv 0$ on the cylinders $T^2\times [-(r-r_0),r-r_0]$, and
the lower and upper bounds by constants $c,C'>0$ independent of 
$r\geq r_0$ follow by the pointwise bound on $\psi_r$ and
$F_{A_r}$. The constant $C'$ depends on the perturbation $\mu$. 

Now consider the case of a compact set
$Z=V\cup_{T^2} [0,r_0]\times T^2$ in $Y(r)$. Applying \eqref{H2} of
the previous Lemma we estimate
$$ c \leq - \frac{\kappa(Y(r_0))}{4} \| \psi_r \|^2_{L^2(Z)}
\leq {\cal E}_{Z,\mu}(A_r,\psi_r) 
 = f|_{\partial Z} + \int_Z F_{A_r}
\wedge \mu, $$
with the boundary term 
$$ f|_{\partial Z}=f(a_r,\alpha_r,\beta_r)= -\int_{\partial Z} \langle
\alpha_r, i\bar\partial^* \beta_r \rangle \, dv_{T^2}. $$
By the assumptions on the metric and on the perturbation we know that
on $Z^c=Y(r)\backslash Z$ we have ${\cal E}_{Z^c,\mu} \geq 0$, and
$\mu\equiv 0$, hence ${\cal E}_{Z^c,\mu}=f|_{\partial
Z^c}=-f|_{\partial Z} \geq 0$. Moreover, for $Z$ of the form as above,
we have 
$$ -C' \leq  \int_Z F_{A_r}
\wedge \mu= \int_{Y(r_0)} F_{A_r}
\wedge \mu \leq C'. $$
Thus we have an estimate
$$ c-C' \leq f(a_r,\alpha_r,\beta_r) \leq 0. $$

In the case of the cylinder region $Z_r$, by considering the two
components in the complement $Z_r^c$ and arguing as above, we obtain
uniform bounds on $f(a_r(r),\alpha_r(r),\beta_r(r))$ and
$f(a_r(-r),\alpha_r(-r),\beta_r(-r))$. The variation
$$ f(a_r(r),\alpha_r(r),\beta_r(r))-
f(a_r(-r),\alpha_r(-r),\beta_r(-r)) ={\cal E}_{Z_r,\p}(A_r,\psi_r) $$
is therefore uniformly bounded in $r\geq r_0$.
\endproof

\begin{Lem} 
Suppose we are given solutions $(A_r,\psi_r)$ of the perturbed 
Seiberg-Witten equations (\ref{3dSW}) 
on the compact 3-manifold $Y(r)$, with a perturbation $\mu$ supported
in the knot complement in $Y(r)$, for $r\geq r_0$. Suppose given a
compact set $Z$ of the form  $V\cup_{T^2} 
[0,r_0]\times T^2$  or $\nu(K) \cup_{T^2}[-r_0,0]\times T^2$ in $Y(r)$
with  $r>r_0$. Then we have uniform bounds 
\[ \| \nabla_{A_r} \psi_r \|^2_{L^2(Z)} \leq C(\kappa,\mu), \ \ \ \ \
\| F_{A_r} \|^2_{L^2(Z)} \leq C(\kappa,\mu) \]  
where $C(\kappa,\mu)$ is a positive constant, depending on the
scalar curvature and on the perturbation, independent of
$r\geq r_0$. 
\label{nablabound}
\end{Lem} 

\proof 
In order to derive the estimate for the $L^2$--norm of $\nabla_A
\psi$, we use the result of Lemma \ref{Hbound}. We have
$$ c\leq \| \nabla_{A_r} \psi_r \|^2_{L^2(Z)} + \int_Z
\frac{\kappa}{4} |\psi_r |^2 dv \leq C' $$
and
$$ c\leq \frac{1}{2} \| F_{A_r} \|^2_{L^2(Z)} + \int_Z
\frac{\kappa}{4} |\psi_r |^2 dv \leq C'. $$
Since the second term is uniformly bounded in $r\geq r_0$, we obtain
the result.
\endproof 

Notice that the uniform bound on the curvature justifies our choice of
the finite energy condition \eqref{finite-energy} for monopoles on the
manifold $V$ with infinite cylindrical end.

Now we can establish the convergence result for the Seiberg-Witten  
monopoles on $Y(r)$ as $r\to \infty$. 
 
\begin{Pro} 
Assume that the metric on $Y(r)$ and the perturbation are 
chosen as specified in the beginning of this section. 
Suppose the moduli spaces ${\cal M}^*_{Y(r)}( \s_r)$ are non-empty for 
$r>>0$, and let $(A_r,\psi_r)$ be a solutions representing elements 
in ${\cal M}^*_{Y(r)}( \s_r)$.    
\begin{enumerate} 
\item 
For any fixed compact set $Z= V\cup_{T^2} (T^2\times
(0,r_0]) \subset Y(r)$,   
there exist gauge transformations $\lambda_r$ on $Y(r)$, such that 
a subsequence of $\lambda_r (A_r,\psi_r)$ converges smoothly 
on $Z$ to either a solution 
$(A',\psi')$ with $[A',\psi']$ in ${\cal 
M}_V ( \s')$, or to a solution $(a_\infty '',0)$, with 
$[a_\infty '',0]$ in  
$ {\cal M}^{red}_{\nu(K)} (\s'')=\chi (\nu(K) , Y ). $ 
\item 
The solutions $\lambda_r(A_r,\psi_r)$ restricted to the cylinder 
$[-r,r]\times T^2$ converge smoothly on compact sets to a constant 
flat connection $a_\infty$ on $T^2$. 
\item 
Let $\partial_\infty [A',\psi']= a_\infty '$ be the asymptotic limit, 
that is, an element of $\chi(T^2,V)$. 
Then there exist two gauge transformations $\lambda'$ and $\lambda''$ on 
$T^2$ that extend to $V$ and $\nu(K)$ respectively, such 
that we have $\lambda'' a_\infty ''= \lambda' a_\infty'$ in 
$\chi(T^2,Y)$. 
\item 
In the case of $Y_0$, we obtain similarly  
$\lambda'' a_\infty ''= \lambda' a_\infty'$ in $\chi(T^2,Y_0)$.  
The gauge transformation 
$(\lambda')^{-1}\lambda''$ over $T^2$ determines a cohomology class in 
$H^2(Y_0, \Z)$ which is the element 
uniquely associated to the $\spinc$ structure $\s$ on $Y_0$. 
\end{enumerate} 
\label{smoothconv} 
\end{Pro} 

\proof

(a) Suppose we are given a fixed compact set $Z= V\cup_{T^2} (T^2\times
(0,r_0])$   
in $Y(r)$. We show that a sequence of elements $[A_r,\psi_r]$ of ${\cal 
M}(Y(r), \s_r)$ 
has a subsequence that converges smoothly on $Z$ to a solution of the 
equations. The same result holds for compact sets $Z$ of the form 
$T^2\times [-r_0, r_0]\cup_{T^2} \nu(K)$. These results were essentially 
established in \cite{KM1}. 
 
The estimates of Lemma \ref{Hbound} and Lemma \ref{nablabound} show
that there is a uniform bound for the $L^2$ norms 
\[ \| \psi_r \|_{L^2 (Z)} \leq C(\kappa, \mu), \ \ \ \ \ \| F_{A_r}
\|_{L^2 (Z)}  \leq C(\kappa, \mu).  \] 
This implies an  $L^2_1$ bound on the connections 
\[ \| A_r-A_0 \|_{L^2_1(Z)}\leq \tilde C\cdot C(\kappa, \mu), \] 
with the constant $\tilde C$ depending on the fixed compact set $Z$,
and independent of $r\geq r_0$. 
The bound 
\[ \| \nabla_{A_r} \psi_r  \|_{L^2 (Z)} \leq C(\kappa,\mu) \] 
of Lemma \ref{nablabound}, together with the $L^2$ bound on the spinors, 
implies a bound on the $L^2_1$-norms of the 
spinors by the elliptic estimate. 

Notice that here $Z$ is a compact
set of the  following form
 \[
V\cup_{T^2} [0,r_0]\times T^2; \qquad  or \qquad \nu(K)
\cup_{T^2}[-r_0,0]\times T^2,
\]
 thus we have elliptic
estimates in the form
$$ \| \psi_r \|_{L^2_{k}(Z)} \leq C ( \| \nabla_{A_r} \psi_r
\|_{L^2_{k-1}(Z')} + \| \psi_r \|_{L^2_{k-1}(Z)} ), $$
where $Z'$ is a smaller set $Z'=V\cup_{T^2} [0,r_0']\times T^2$, for
some $r_0' < r_0$, cf.~\cite{MW} \S 4.1. 
Since we are only taking estimates on a fixed compact set
$Z$ of the form specified above, the constant $C$ in the elliptic
estimate depends on $Z$ but does not depend on the parameter $r$ of
the underlying manifold $Y(r)$. For the elliptic estimate for the
connections, we choose any smooth connection $A_0$ on $det(\s)$ over $Z$ 
and gauge transformations $\lambda_r$ in the identity connected component 
of the gauge group ${\cal G}_{Y(r)}$, such that the forms 
$\lambda_r A_r -A_0$ are co-closed and annihilate 
the normal vector at the boundary $T^2$. We use an elliptic
estimate of the form considered above for the operator $d +
d^*$. Thus, we can bound the $L^2_2$-norms of the connections on $Z$,
and use a bootstrapping argument to bound the higher Sobolev norms as
in \cite{KM1} \cite{MW} \cite{Mor}.  

Upon passing to a subsequence, we have obtained elements 
$(A_{r_i},\psi_{r_i})$ that converge smoothly on $Z$ to a solution of 
the equations. This defines a solution $(A',\psi')$ on $V$ 
with the cylindrical end $T^2\times [0,\infty)$. 
The case of $\nu(K)$ is analogous. With our choice of metric on $\nu(K)$,  
a finite energy solution on $\nu(K)$ will necessarily be reducible. 

To complete the proof of (a) we need to show that the resulting
solution on $V$ with infinite end satisfies the finite energy
condition \eqref{finite-energy}. This follows from the curvature estimate
in Lemma \ref{nablabound}

(b) To prove the second claim, consider the elements $x_r=(A_r,\psi_r)$ 
restricted to the cylinder  
$[-r,r]\times T^2$. Up to a gauge transformation, they can be written 
in the form 
$$ x_r(s)= (a_r(s),\alpha_r(s),\beta_r(s)).$$ 

The functional $f$ is monotone along the cylinder, with variation
$$
f(a_r(r),\alpha_r(r),\beta_r(r))-f(a_r(-r),\alpha_r(-r),\beta_r(-r))=
$$
$$\int_{-r}^r \| \nabla f(a_r(s),\alpha_r(s),\beta_r(s)) \|^2 ds. $$
By the result of Lemma \ref{Hbound}, there is a uniform bound,
independent of $r$ for the variation  of the  
functional $f$ along the cylinder, 
\[  f(a_r(r),\alpha_r(r),\beta_r(r))- 
f(a_r(-r),\alpha_r(-r),\beta_r(-r)) \leq C. \] 
%In fact, when considering the manifold $Y(r)=V\cup_{T^2} [-r,r]\times
%T^2 \cup_{T^2} \nu(K)$, we identify the section $T^2\times \{ -r \}$
%of the cylinder in $Y(r)$ with $T^2\times \{ 0 \}$ of $V(r)\subset
%Y(r)$, and the section $T^2\times \{ r \}$ in $Y(r)$ with $T^2\times
%\{ 0 \}$ of $\nu(K)(r)\subset Y(r)$. By the smooth convergence on
%compact sets of part (a), applied to $V=V(0)$ and $\nu(K)=\nu(K)(0)$,
%we obtain the uniform bound on the values of $f(a_r,\alpha_r,\beta_r)$
%on $\partial V$ and $\partial \nu(K)$. 

This uniform bound for 
$$ c\leq {\cal E}_{Z_r,\mu}(A_r,\psi_r) \leq C, $$
with $Z_r=T^2\times [-r,r]$, gives bounds on compact sets
$Z\subset Z_r$ for the $L^2$--norms 
$\| \nabla_{A_r}\psi_r \|$, $\| \psi_r \|$, $\| F_{A_r} \|$, as well as 
for the $L^4$ norm of the spinor. This is enough to start the
bootstrapping argument, with elliptic estimates as before, hence we
obtain smooth convergence on compact sets in $Z_r$ to a solution of
the unperturbed SW equations on $T^2\times \R$. Such solution must be
a flat connection and the trivial spinor. This implies ${\cal
E}_{Z_r,\mu}(A_r,\psi_r) \to 0$, hence, using again Lemma
\ref{Hmonopoles} together with the estimate (\ref{Loj-f}), we obtain
that the limit is actually a critical point $a_\infty$ of $f$. 

Thus, up to gauge transformations, the sequence of 
solutions $(A_r,\psi_r)$ has a subsequence $(A_{r_i},\psi_{r_i})$ that 
converges smoothly on compact sets to a pair 
$((A',\psi'),(a_\infty '',0))$. In the asymptotic limit we get 
$$\lim_{s\to\infty} \lambda' (A',\psi')= \lambda' a_\infty'=  
\lambda'' a_\infty''.$$ 
 
In the case of the manifold $Y_0$, 
\[
[\lambda' a_\infty']= 
[\lambda'' a_\infty'']
\]
 in $\chi(T^2, Y_0)$ imply that 
$x_r= [A_r, \psi_r] \in \M^*_{Y_0(r)}(\s_k)$ where $\s_k$ corresponds to the
cohomology class 
\[ 
[(\lambda')^{-1} \lambda''] \in H^1(T^2, \Z)/ H \]
for 
\[ H = Im(H^1(V, \Z), H^1(T^2, \Z)) + Im(H^1(\nu(K), \Z), H^1(T^2,
\Z)). \] 
\endproof

This completes the proof of the convergence part of the gluing 
theorem \ref{Gluing} for generators. Namely, we have shown
that a gauge class in the moduli space  ${\cal M}_{Y(r)}(\s)$, for a
sufficiently large $r$, and perturbation as prescribed, determines an
element in  
\[ {\cal M}^*_{V, Y} \times_{\chi(T^2, Y)} \chi (\nu(K), Y). \]

\subsection{Proof of Theorem \ref{Gluing}}

In this subsection we will construct an approximate monopole on 
$Y(r)$ from any element  
in $ {\cal M}^*_{V, Y} 
\times_{\chi(T^2, Y)} \chi (\nu(K), Y)$, 
and study the gluing that will produce the 
corresponding monopole on $Y(r)$ for a  sufficiently large $r$.  
 
First, we define a pre-gluing operation, where we splice together solutions in 
${\cal M}^*_V$ and $\chi (\nu(K))$ with matching asymptotic 
values, via a smooth cutoff function. This produces an approximate 
solution $(A',\psi')\#_r (a_\infty'',0)$ of the monopole equations on $Y(r)$ 
for $((A',\psi'), (a_\infty'',0))$ representing an element in 
$ {\cal M}^*_{V, Y} 
\times_{\chi(T^2, Y)} \chi (\nu(K), Y)$.

We can assume that $(A',\psi')|_{T^2\times [0, \infty) \subset V}$ 
is in temporal gauge with asymptotic limit $(a_{\infty}, 0)$, and there 
is a gauge transformation $\lambda''$ on $\nu(K)$ such that  
$\lambda'' (a_\infty'') = a_{\infty}$ as a flat connection on $T^2$. 
Let $(A',\psi') = a_\infty + (a'(s), \psi'(s))$ on $T^2\times [0, \infty)$. 
We can choose smooth cutoff functions $\rho_r  (s)$ ($s\in [-2, 2]$) 
with values in $[0,1]$, satisfying $\rho_r (s)\equiv 1$ for  
$s\in [-2, -1]$  and $\rho_r (s)\equiv 0$ for $s\in [1, 2]$ 
with $0\leq \rho'(s) \leq 1$. 
 
Define the pre-gluing map with values in ${\cal B}(Y(r))$ 
by setting 
\ba  \begin{array}{lll}
&&x_r = (A, \psi) \\[2mm]
&=&(A',\psi')\#^0_r (a_\infty '',0)\\[2mm]
& =&  
\left\{ \begin{array}{lr} 
(A',\psi') & \hbox{on}~ V(r-2) \\ 
a_\infty +\rho_r  (s)\lambda''(a'(s+r),\psi'(s+r))  
&  s\in [-2, 2] \\ 
\lambda'' (a_\infty '',0) & on~ \nu(K)(-r +2) 
\end{array}\right. \end{array} 
\label{pre-map}  
\na 
 
\begin{Def} 
An {\bf approximate solution} is by definition an element in the image of 
the pre-gluing map (\ref{pre-map}).  
We use the notation
\[
{\cal M}^*_{V, Y} 
(a_\infty):=\partial_\infty^{-1}(a_\infty) \subset {\cal M}^*_{V, Y}. 
\]
 Then  ${\cal U}(a_\infty, r)$ is defined to be 
 the image of the pregluing map 
(\ref{pre-map})  
$ \#^0_r: {\cal M}^*_{V, Y}(a_\infty)\times [a_\infty '', 0] \to {\cal 
B}(Y(r)).$ 
\label{approximate-def} 
\end{Def}

In order to show that the approximate solutions in ${\cal U}(a_\infty, 
r)$ can be deformed to actual solutions of the monopole equations on 
$Y(r)$, we consider the span of 
eigenvectors corresponding to the small eigenvalues of  
the linearization operator at the approximate solutions.  
 
Consider the linearization operator of the Seiberg-Witten equations 
on $Y(r)$ at the approximate solution $(A',\psi')\#^0_r 
(a_\infty '',0)$ 
$$ H_{(A',\psi')\#^0_r  
(a_\infty '',0)}(f,\alpha,\phi)= \left\{\begin{array}{l} 
            L_{(A',\psi')\#^0_r 
(a_\infty '',0)}(\alpha,\phi)+G_{(A',\psi')\#^0_r (a_\infty '',0)}(f)\\ 
            G^*_{(A',\psi')\#^0_r (a_\infty '',0)}(\alpha,\phi). 
          \end{array}\right.  $$ 
We also need the linearization operators of the Seiberg-Witten equations 
on $V$ and $\nu(K)$ with infinite cylindrical ends, as defined in the
deformation complex (\ref{deform-fiber}), acting on $L^2$ forms and
spinors:  
$$ H_{(A',\psi')}(f,\alpha,\phi)= \left\{\begin{array}{l} 
            L_{(A',\psi')}(\alpha,\phi)+G_{(A',\psi')}(f)\\ 
            G^*_{(A',\psi')}(\alpha,\phi) 
          \end{array}\right. $$ 
$$ H_{(a_\infty '',0)}(f,\alpha,\phi)= \left\{\begin{array}{l} 
            L_{(a_\infty '',0)}(\alpha,\phi)+G_{(a_\infty '',0)}(f)\\ 
            G^*_{(a_\infty '',0)}(\alpha,\phi) 
          \end{array}\right. , $$ 
where the operator $L$ is defined as in (\ref{L}). We think of 
$H_{(A',\psi')\#^0_r (a_\infty '',0) }$ as acting on the elements 
$(\alpha,\phi)$ in the $L^2$ tangent space of the configuration space 
over the closed manifold $Y(r)$.  
We continue to  denote by $H_{(A',\psi')}$ and $H_{(a_\infty '',0)}$ the 
operators defined in the deformation complex 
(\ref{deform-M}) acting on the extended $L^2$ spaces of connections 
and spinors, over $V$ and $\nu(K)$ respectively. 
 
Now we discuss the  
eigenfunctions corresponding to slowly decaying eigenvalues 
of these operators. 
The model for our analysis of the operator 
$H_{(A',\psi')\#^0_r(a_\infty'',0)}$ is based on the work of Capell, 
Lee, and Miller \cite{CLM1},\cite{CLM2}. 
With operators differing from a translation invariant
operator by exponentially decaying terms, we shall adopt the more
general setting as in the work of Nicolaescu, \cite{Nicol}.   
 
We use the following result, which is the analogy in our context  
of Theorem A of \cite{CLM1}.  
 
\begin{Pro} 
Assume that $a_\infty$ is a point in $\chi(T^2)$ away from a small 
neighbourhood $U_\Theta$ of $\Theta$. 
Let   
\[
N(r)= \dim Ker_{L^2}(H_{(A',\psi')}) + \dim Ker_{L^2}( 
H_{(a_\infty '',0)}) + \dim Ker(Q_{a_\infty}).
\]
Then, there exists  an $N(r)$-dimensional family of eigenvectors of 
the operator $ H_{(A',\psi')\#^0_r (A'',\psi'')}$ with eigenvalues 
satisfying $\lambda(r)\to 0$ as $r\to\infty$ 
at the rate at most $1/r$. 
The dimension $N(r,r^{-(1+\epsilon)})$ of the span of eigenvectors of 
the operator $ H_{(A',\psi')\#^0_r (a_\infty '',0)}$ with eigenvalues 
$\lambda < r^{-(1+\epsilon)}$ is given by 
$$ N(r,r^{-(1+\epsilon)})= \dim Ker_{L^2}(H_{(A',\psi')}) + \dim Ker_{L^2}( 
H_{(a_\infty '',0)}) + \dim \ell_1 \cap \ell_2, $$ 
where $\ell_1$ and $  \ell_2$ are the two Lagrangian 
submanifolds in $Ker(Q_{a_\infty}) = H^1(T^2,\R)$, determined by the 
extended $L^2$ solutions of $H_{(A',\psi')}(\alpha,\phi)=0$ and 
$H_{(a_\infty '',0)}(\alpha,\phi)=0$. 
\label{CLM} 
\end{Pro} 
 
\begin{proof} 
In order to prove the first claim it is sufficient to 
check that elements of $H^1(T^2,\R)=Ker(Q_{a_\infty})$ give rise to 
approximate eigenfunctions on $Y(r)$ with slowly decaying 
eigenvalues, that is, with eigenvalues  $\lambda(r)$ satisfying 
$\lambda(r)\to 0$ at most like $1/r$. The first statement is then an 
analogue, in our case, of Proposition 6.B of \cite{CLM1}. 
 
Suppose we are given an element $\xi\in Ker(Q_{a_\infty})$. 
If $\chi(s)$ is a cutoff function supported 
in $[r/2-\epsilon, 3r/2+\epsilon]$ satisfying $\chi(s)\equiv 1$ on 
$[r/2,r]$, we have an estimate 
\[ \frac{\| (\d_s + Q_{a_\infty}) \chi \xi \|_{L^2(Y(r))}}{\|\chi \xi 
\|_{L^2(Y(r))}} \leq \frac{C}{r}. \] 
This implies a similar estimate for the operator $H_{(A',\psi')\#^0_r 
(a_\infty '',0)}$ on $Y(r)$, for $r\geq r_0$ large enough, since we 
are assuming that this operator differs from $\partial_s + 
Q_{a_\infty}$ by terms that are exponentially small in $r$. 
This is the setting used in \cite{Nicol}.  
 
The second part of the statement 
can be derived from the asymptotically exact sequence 
$$ 0\to {\cal K}(r^{-(1+\epsilon)}) \to 
Ker_{L^2}^{ext}(H_{(A',\psi')}) \oplus 
Ker_{L^2}^{ext}(H_{(a_\infty'',0)}) \stackrel{\Delta}{\to} \ell_1 
\oplus \ell_2 \to 0, $$ 
as in the Main Theorem of \cite{Nicol}.  
Here ${\cal K}(r^{-(1+\epsilon)})$ denotes the span of the eigenvectors of  
$ H_{(A',\psi')\#^0_r (A'',\psi'')}$ with small eigenvalues that decay 
at a rate of at least $r^{-(1+\epsilon)}$. We use the notation 
$Ker_{L^2}^{ext}$ for the extended $L^2$-solutions, and $\ell_i$ for 
the asymptotic values of the extended $L^2$-solutions. 
\end{proof} 
 
Proposition \ref{CLM} yields the following. 
 
\begin{Cor} 
There are no fast decaying eigenvalues, that is, in our problem 
$N(r,r^{-(1+\epsilon)})=0$. However, there is a non--trivial family of
eigenvectors of the linearization  
$ H_{(A',\psi')\#^0_r (a_\infty '',0)}$ at the approximate solution 
$(A',\psi')\#^0_r (a_\infty '',0)$, with slowly decaying 
eigenvalues, satisfying  
$\lambda(r)\to 0$ at most like $1/r$. 
\label{low-modes} 
\end{Cor} 
 
\begin{proof} 
We have 
$$ \dim Ker_{L^2}(H_{(A',\psi')})=\dim Ker_{L^2}(H_{(a_\infty '',0)})=0. $$ 
Moreover, for a generic choice of the perturbation of the 
monopole equations on $V$, the Lagrangian subspaces $\ell_1$ and 
$\ell_2$ intersect transversely. 
Thus, we have $N(r)=\dim Ker(Q_{a_\infty})$ and 
$N(r,r^{-(1+\epsilon)})=0$. The previous Proposition shows that the
span of eigenvectors with slowly decaying eigenvalues is
non--trivial. In fact, it shows the existence of (at least) a two
dimensional family parameterized by the elements of
$H^1(T^2,\R)=Ker(Q_{a_\infty})$. 
 
\end{proof}

Suppose we are given an element $(a, \phi)$ on $Y(r)$ such that 
$x_r +(a, \phi)$ is a solution of the monopole equations on 
$Y(r)$. Then $(a, \phi)$ satisfies 
$$ H_{x_r}(a,\phi)+ N_{x_r} (a,\phi) + 
\Sigma(x_r) =0, $$ 
where $\Sigma$ is the error term defined by
$$ \Sigma(x_r)=\left( \begin{array}{c} *F_A-\sigma(\psi,\psi)-\mu \\ 
\dirac_A \psi \end{array}\right), $$ 
as by equation \eqref{3dSW}, and $N$ is the non--linear term
$$ N_{A,\psi}(a,\phi)=\left( \begin{array}{c} \sigma(\phi,\phi) 
\\ a .\phi \end{array}\right). $$

Though we do not treat the more general case in this paper, we 
mention that one can consider the same argument with an additional
perturbation term $P(A,\psi)$. In this case, an additional term
$P(A,\psi)$ enters the expression for the error term $\Sigma(x_r)$,
and an additional non-linear part ${\cal N}P_{A,\psi}$ of the
perturbation 
$$ {\cal N}P_{A,\psi}=P((A,\psi)+(a,\phi))-{\cal 
D}P_{A,\psi}(a,\phi) $$ 
is added to the expression of $N_{A,\psi}$. This case will be
discussed elsewhere.

Choose $\lambda=\lambda(r)>0$ such that $\lambda(r)$ is not an eigenvalue of  
$H_{x_r}= H_{(A',\psi')\#^0_r (a_\infty '',0)}$, for all approximate
solutions $x_r =(A',\psi')\#^0_r (a_\infty '',0)$ in ${\cal U}(a_\infty, 
r)$. Consider the projection maps $\Pi(\lambda(r),x_r)$ onto the span of 
the eigenvectors of $H_{x_r}$ with eigenvalues smaller than $\lambda(r)$. 

The condition that, for a given approximate solution $x_r$, the
element $x_r +(a, \phi)$ is an actual solution of monopole 
equations can be written as a system of two equations: 
\ba 
\Pi(\lambda(r),x_r)(N(a, \phi) + 
\Sigma(x_r)) =0  
\label{fix-1} 
\na  
\ba 
H_{x_r}(a, \phi)+ 
(1-\Pi(\lambda(r),x_r))(N(a, \phi) + \Sigma(x_r)) =0.  
\label{fix-2} 
\na

If the equation (\ref{fix-2}) admits a unique solution 
$(a, \phi)$, then the  
condition that $x_r +(\alpha, \phi)$ is a solution of the monopole 
equations on $Y(r)$ can be rephrased as the condition that 
(\ref{fix-1}) is satisfied, with $(a, \phi)$ the unique solution 
of (\ref{fix-2}). 

The second equation \eqref{fix-2} can be written as the fixed point
problem  
\begin{equation} \label{fixedpoint} (a, \phi)= - H_{x_r}^{-1} 
(1-\Pi(\lambda(r),x_r))(N(a, \phi) +  \Sigma(x_r)). \end{equation} 
The following result proves existence and uniqueness of the solution
to \eqref{fixedpoint}. 
 
\begin{Lem} 
There exists a positive constant $C >0$, such that,
if a given approximate solution $x_r$ satisfies
$\| \Sigma(x_r) \|_{L^2(Y(r))} \leq C \epsilon(r)^2$, 
for some small and positive
$\epsilon(r)$ satisfying $ \epsilon(r)< \frac{\lambda(r) }{2C}$, then
the map
$$ {\cal T}_r (a, \phi):= - H_{x_r}^{-1} 
(1-\Pi(\lambda(r),x_r))(N(a, \phi) +  \Sigma(x_r)) $$ 
maps the ball $B_{\epsilon(r)}=\{ (a, \phi) | \ \| (a, \phi) 
\|^2_{L^2_1(Y(r))} \leq \epsilon(r) \}$ to itself and is a contraction
on $B_{\epsilon(r)}$. 
\label{contraction} 
\end{Lem} 

\begin{proof} Let $C>0$ be a constant such that the quadratic term
satisfies the estimate
$$ \| N(a,\phi)- N(\tilde a,\tilde\phi) \|_{L^2} \leq C ( 
\|(a,\phi) \|_{L^2_1} + \| (\tilde a,\tilde\phi) \|_{L^2_1}) \| 
(a,\phi)- (\tilde a,\tilde\phi) \|_{L^2_1}, $$   
independent of $r\geq r_0$. This follows from the Sobolev multiplication
theorem in dimension 3.

On the image of $(1-\Pi(\lambda(r),x_r))$, the operator 
$H_{x_r}^{-1}$ is bounded with  norm bounded by $\lambda(r)^{-1}$. 
We have an estimate for $(a, \phi)\in B_{\epsilon(r)}$
\[\begin{array}{lll}
 && \| {\cal T}_r (a, \phi) \|_{L^2_1(Y(r))}\\[2mm]
&\leq& \frac{1}{\lambda(r)} \| N (a, \phi) + \Sigma (x_r)\|_{L^2(Y(r))}
\\[2mm] &\leq& \frac{C\epsilon(r)^2}{\lambda(r)}  + 
\frac{\|\Sigma (x_r)\|_{L^2(Y(r))}}{\lambda(r)} \\[2mm] 
&\leq& \frac{2C\epsilon(r)^2}{\lambda(r)} \leq \epsilon(r),
\end{array}\]
which implies that ${\cal T}_r$ maps the ball $ B_{\epsilon(r)}$ to itself.  

Let $(a_1, \phi_1), (a_2, \phi_2) \in B_{\epsilon(r)}$, we have 
\[
\begin{array}{lll}
&& \| {\cal T}_r (a_1, \phi_1)- {\cal T}_r (a_2, \phi_2) \|_{L^2_1(Y(r))}
\\[2mm]  
&\leq& \frac{1}{\lambda(r)} \| N (a_1, \phi_1)- N (a_2, \phi_2) \|_{L^2} 
\\[2mm]
&\leq& \frac{C}{\lambda(r)} \| (a_1, \phi_1) + (a_2, \phi_2)\|_{L^2_1}
\| (a_1, \phi_1) -(a_2, \phi_2)\|_{L^2_1}\\[2mm]
&\leq& \frac{2C\epsilon(r)}{\lambda(r)}
\| (a_1, \phi_1) -(a_2, \phi_2)\|_{L^2_1}. 
\end{array}
\] 
Thus, from $\epsilon(r) < \frac{\lambda(r)}{2C}$ as chosen,  we obtain 
 that ${\cal T}_r$ is a contraction on $B_{\epsilon(r)}$.
\end{proof}

\begin{Pro}
For sufficiently large $r\geq r_0$, and for all approximate solutions
$x_r$ in ${\cal U}(a_\infty, r)$, there exists a unique solution $(a,
\phi)$ of \eqref{fix-2}, such that equation \eqref{fix-1} is trivially
satisfied. 
\end{Pro} 

\begin{proof} For all approximate solutions $x_r$ in ${\cal
U}(a_\infty, r)$, we have an estimate on the error term
$$ \| \Sigma(x_r) \|_{L^2(Y(r))} \leq C' e^{-\delta r}, $$ 
for $r\geq r_0$, which follows from the exponential decay estimate
proved in Proposition \ref{decay-f}.  
Thus, we can apply Lemma \ref{contraction}, with $\lambda(r) =
O(r^{-(1+\epsilon)})$ and $\epsilon(r)= O(e^{-\delta r/2})$.
By Corollary \ref{low-modes} we know that, for $\lambda(r) =
O(r^{-(1+\epsilon)})$, the projection $\Pi(\lambda(r),x_r)\equiv 0$, hence
the solution $(a,\phi)$ of \eqref{fix-2}, provided by Lemma
\ref{contraction} also satisfies trivially equation \eqref{fix-1}.

\end{proof} 

Thus, the resulting element $x_r + (a,\phi)$ is a true monopole
solution on $Y(r)$, close to the approximate solution $x_r$. This
completes the proof of the gluing theorem \ref{Gluing}.

\subsection{Metric}
 
Concerning the metric after surgery,  
on $\nu(K)$ inside $Y_1$ we consider the metric 
$g_1$ as in the following Lemma \ref{Liviu-metrics}, 
which is due to Liviu Nicolaescu \cite{Ni-pri}. 
 
Let $g = du^2 + dv^2$ with $\int_{T^2} du\wedge dv = 4\pi^2$, where
the torus $T^2$ is the  
boundary of the tubular neighbourhood of the 
knot $\nu(K)$ in $Y$. We introduce a 
choice of a metric on $\nu(K)$ inside $Y_1$, for which we can still 
derive the result that the moduli space of monopoles on $\nu(K)$ inside $Y_1$ 
consists of the circle of reducibles.  
 
\begin{Lem} (Nicolaescu)
Let $A$ be an element in $SL(2,\Z)$. Suppose we are given $\epsilon >0$ 
sufficiently small. Consider the flat metric on $T^2$ given by 
$ g_0= A^* g$. There exists a 
constant $c$ and a smooth path $g(s)$ ($s\in \R$) of flat metrics on $T^2$ 
with the following properties: 
 
(i) $g(s)\equiv \frac{1}{\delta^2} g_0$, for all $s \leq \epsilon$ 
and $g(s)=g_1$ for all $s\geq 1-\epsilon$; 
 
(ii) $g_1 = g(1)$ is a metric of the form 
$ g_1 = k_1 du^2 + k_2 dv^2 $ 
with positive constants $k_i$; 
 
(iii) The scalar curvature of the metric $\hat g:= g(s)+ ds^2 $ on $T^2 
\times \R$ is non-negative; 
 
(iv) The metric $g_1$ can be extended to a metric inside the 
solid torus $\nu(K)$ with positive scalar curvature. 
\label{Liviu-metrics} 
\end{Lem} 

\noindent \underline{Proof of Nicolaescu's Lemma}:
Choose a unit vector $\partial_u$ with respect to the metric 
$\frac{1}{c} g_0$, and complete it to an oriented orthonormal 
frame. Let $\{ \varphi_1, \varphi_2 \}\subset \Omega^1(T^2)$ be the 
dual coframe. This is related to $\{ du, dv \}$ by 
$$ \varphi_1 = du + a_0 dv \ \  \varphi_2=k \, dv, $$ 
for some positive constant $k>0$. 
 
The path $g(s)$ is defined by requiring that the coframe 
$$ \varphi_1(s) = du + a(s) dv \ \  \varphi_2(s) =k \, dv $$ 
be orthonormal with respect to $g(s)$, where $a(s)$ is a smooth 
function satisfying $a(s)\equiv 0$ for all $s\geq 1-\epsilon$ and 
$a(s)=a_0$ for all $s\leq \epsilon$. The only conditions that need to 
be verified are (iii) and (iv). 
 
We have an orthonormal coframe $\{ \varphi_0, \varphi_1, \varphi_2 \}$ 
on $X= T^2\times\R$, with respect to the metric 
$\hat g$, with $\varphi_0=ds$, and a corresponding orthonormal frame 
$\{ e_0, e_1, e_2 \}$. The Levi-Civita connection is of the form 
$$ \Gamma=\left( \begin{array}{ccc} 0&x&y\\ -x&0&z\\ -y&-z&0 
\end{array}\right) \ \ x,y,z\in \Omega^1(X). $$ 
 
The Cartan structural equation gives 
$ d \vec{\varphi} =\Gamma \wedge\vec{\varphi}$, 
with $\vec{\varphi}=(\varphi_0, \varphi_1, \varphi_2)$. By the 
expression of $\varphi_i$, we have 
$$ d\varphi_0=d\varphi_2=0, \ \ d\varphi_1= \frac{\dot{a}}{k} 
\varphi_0\wedge \varphi_2, $$ 
hence we obtain 
\[\begin{array}{lll}  
&& x\wedge \varphi_1+ y\wedge \varphi_2=0, \\[2mm] 
 &&\frac{\dot{a}}{k} \varphi_0\wedge \varphi_2= -x \wedge\varphi_0 + 
z\wedge \varphi_2, \\[2mm] 
&&  -y \wedge \varphi_0 - z\wedge \varphi_1=0.\end{array} 
\] 
These equations imply 
 $$ \Gamma= \frac{\dot{a}}{2k}\left( \begin{array}{ccc} 0&\varphi_2 & 
\varphi_1 \\ -\varphi_2 & 0 & \varphi_0 \\ -\varphi_1 & -\varphi_0 & 0 
\end{array} \right). $$ 
Thus, we can compute the scalar curvature of $\hat g = g(s) + ds^2$  
on $T^2 \times \R$ which is  
$3(\dot{a}/k)^2$ by direct calculation. 
 
Claim (iv) then follows by noticing that any diagonal metric of the form 
$$ g_1 = k_1 du^2 + k_2 dv^2 $$ 
realizes the torus $T^2$ metrically as the product of two circles of 
different radii. Each can bound a hemisphere, endowed with a positive scalar 
curvature metric, thus extending $g_1$ to a metric on a solid torus, 
with positive scalar curvature. 
 
\endproof

\subsection{Lines in $\chi(T^2)$}

In this subsection we justify why it is sufficient to consider the gluing
map in Theorem \ref{Gluing} away from the singular point $\Theta\in
\chi(T^2)$. 

\begin{Lem}
The intersection points $\partial_\infty \M^*_V \cap \chi(\nu(K))$,
with $\chi(\nu(K))\subset \chi(T^2)$ the circle of reducibles for
$\nu(K)$ in either $Y$, $Y_1$, or $Y_0$, are contained in
$\chi(T^2)\backslash U_\Theta$, for some neighborhood $U_\Theta$ of
the singular point $\Theta$.
Thus, the gluing of Theorem \ref{Gluing} happens away from the
reducible point.
\label{noredgluing}
\end{Lem}
 
\proof 
The torus $T^2$ inside $Y$ inherits from the trivial $Spin$ 
structure of $Y$ the non-trivial $Spin$ structure in which both 
circles (longitude and meridian) bound, that is, the one determined by 
the element $(1,1)$ in $H^1(T^2,\Z_2)$. Introduce the coordinates  
$u$ and $v$ on $H^1(T^2,\R)$,  defined 
by the property that, under the projection to $\chi(T^2)$, they 
satisfy the following condition. 
For $[A]  \in \chi (T^2)$, $v ([A]) $ is the holonomy 
around the meridian $m$ and $u ([A]) $ is the holonomy 
around the longitude $l$. Under this coordinate 
system, the singular point $\Theta$ is given by 
$(1, 1)$ in $\chi(T^2)$, and the reducible circle 
$\chi (V) = \M_V^{red}$ is given by $\{ u=0\}.$ 
Also for the unperturbed Seiberg-Witten equations on $\nu(K) \subset Y$, 
with the metric of non-negative scalar curvature, 
strictly positive away from the boundary, the reducible circle 
$\chi (\nu(K), Y)$ is given by 
\[ 
L_Y := \{ v = 0 \}\subset \chi _0 (T^2, Y). 
\] 
Similarly, in the case of $Y_1$, choose a metric with a long cylinder 
$[-r,r]\times T^2$, which agrees with the original metric on $Y$ when 
restricted to the knot complement $V$, and such that the induced 
metric in the torus neighbourhood $\nu(K)$ is as described in Lemma 
\ref{Liviu-metrics}, then the reducible circle 
$\chi (\nu(K) , Y_1)$  is given by 
\[ 
L_{Y_1} := \{ v-u =1\}\subset\chi _0 (T^2, Y_1). 
\] 
We might have expected that the reducible case would be 
$L_{Y_1} := \{ v-u =0\}$ from the $+1$-surgery manifold 
$Y_1$. 
The shift is due to the fact that $+1$-surgery changes 
the underlying $Spin$ structure $(1, 1)$ by tensoring with a flat
$\Z_2$-bundle of class  
$(0,1)$ in $H^1(T^2,\Z_2)$. In the case of $Y_0$, 
the  reducible circle 
$\chi (\nu(K) , Y_0)$ for the unperturbed equations 
is mapped to a circle $\{ u=0\} \subset \chi(T^2)$.  
This is because 
\[ 
\chi (T^2, \nu(K)) = \chi (T^2, V) = 
\chi(T^2,Y_0) \cong \R \times S^1, 
\] 
so that $\chi (\nu(K)   , Y_0)$ consist of a $\Z$-family of circles 
given by $\{u=2k, k\in \Z\}$ where $u$ is the coordinate of $\R$ 
in $\R \times S^1$.  
The gluing map on the fiber product 
\[ 
{\cal M}^*_{V, Y_0} \times_{\chi(T^2,Y_0)}\{ u=2k \} 
\] 
would correspond to the moduli space $\M_{Y_0}(\s_k)$ where $\s_k$ 
is the $\spinc$ structure with  
$c_1(det \s_k) \in H^2(Y_0, \Z)$. For the trivial $\spinc$ 
structure $\s_0$, 
there would be a circle of reducible monopoles 
in $\M_{Y_0} (\s_0)$ resulting from gluing  
the reducibles $\chi(V, Y) \#\{ u=0\}$,  
we need to introduce a small  perturbation 
inside $\nu(K)$ such that $\chi (\nu(K) , Y_0) = \{ u= \eta\} $ 
where $\eta$ is  small number in $\R$ for the trivial $\spinc$ 
structure $\s_0$. 
 
Clearly, in all the cases, $\chi (\nu(K))$ does not go through 
the singular points $\{\pi^{-1}(\Theta) \}$, hence there is no need to 
consider the gluing map at the singular point $\Theta$. Note that, 
after perturbing the metric inside a compact set on the manifold $V$
with an infinite cylindrical end, we can 
make the open ends in $\M_V$ not limit to any intersection 
points of $\chi(V)=\M^{red}_V =\{ u=0\}$ with any circle $\chi (\nu(K))$ for 
either $Y, Y_1$ or $Y_0$.    

\endproof

Thus, with our choice of metric as in the previous subsection, and 
with the choice of perturbation as in Theorem 
\ref{smooth-dim1},  we see that 
the gluing result stated in Theorem \ref{Gluing} 
holds for the manifolds $Y$, $Y_1$, and $Y_0$. This completes 
the proof of Theorem \ref{Gluing}. 
 
\begin{Rem} Note that gluing the reducible monopoles 
on $V$ and $\nu(K)$ with matching boundary condition 
just gives the extension of the flat connections to the whole 
manifold (after 
a possible gauge transformation). 
We call this the trivial gluing.  
The unique reducible point $\theta_Y$ in ${\cal M}_{Y}$ 
is obtained by the trivial gluing of the 
unique 
intersection point between the lines 
$ L_Y = \{ v = 0 \} $ and  
$ \{ u = 0 \}= \pi^{-1}(\chi(V))\subset \chi(T^2,Y)$.  
The unique reducible point $\theta_{Y_1}$  in ${\cal M}_{Y_1}$ 
is obtained by the trivial gluing of the 
unique 
intersection point between the lines 
$ L_Y = \{ v - u = 1 \}$ and $ \{ u = 0 \}=  
\pi^{-1}(\chi(V))\subset \chi(T^2,Y_0)$. 
\end{Rem} 
 
\section{The geometric triangle and proof of Theorem \ref{surgery-SW3}} 
 
\begin{figure}[ht] 
\epsfig{file=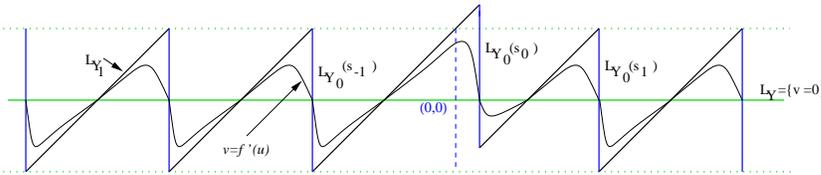,angle=0} 
\caption{The geometric triangles \label{fig0-tri}} 
\end{figure}

In the previous section, we showed that the moduli spaces 
for irreducible monopoles on $Y, Y_1$ and $Y_0$ are given by 
the gluing maps on the following fiber products: 
\ba\begin{array}{l} 
 {\cal M}^*_{Y(r)} \cong {\cal M}^*_{V,Y} 
\times_{\chi(T^2,Y)} \{ v=0\}, \\ 
{\cal M}^*_{Y_1(r)} \cong {\cal M}^*_{V,Y_1} 
\times_{\chi(T^2,Y_1)} \{ v-u =1\}, \\ 
{\cal M}^*_{Y_0(r)} ({\s_k})  \cong {\cal M}^*_{V,Y_0} 
\times_{\chi(T^2,Y_0)} \{ u= 2k\}, \qquad  \qquad \ for\ k\neq 0,\\ 
{\cal M}^*_{Y_0(r)} ({\s_0 })  \cong {\cal M}^*_{V,Y_0} 
\times_{\chi(T^2,Y_0)} \{ u= \eta\}, 
\end{array} 
\label{3moduli} 
\na 
where $(u, v)$ is the coordinate system on $\chi (T^2)$ and 
its covering spaces, $\eta > 0$ is a small parameter, and 
$r>>0$ is a sufficiently large number. We can study these moduli spaces 
on the common character variety $\chi(T^2,Y_0)$ which 
can be identified as a cylinder $\R^1\times S^1$. 
Specifically we take it to be the 
domain (see Figure \ref{fig0-tri}) 
\[ 
\{(u, v)| u\in \R, v \in [-1, 1]\} /\{ (u, -1)\thicksim (u, 1)\} 
\] 
in which the lines corresponding to $L_Y, L_{Y_1}$ and 
$L_{Y_0}(\s_k) (k\in \Z)$ are drawn. 
 
In this section, we introduce a suitable perturbation of the curvature 
equation, supported in the solid torus $D^2 \times S^1$, that 
simulates the effect of surgery such that the reducible line 
corresponding to $\nu(K)\subset Y$ is given by the curve 
$v= f'(u)$ as shown in Figure \ref{fig0-tri}.  
 
For a generic 
perturbation we can assume the curves $\d_{\infty} (M_V^*)$ stay away from 
the intersection  points $\{L_Y\cap L_{Y_1}, L_Y \cap L_{Y_0}, 
L_{Y_1} \cap L_{Y_0}\}$, hence $\d_{\infty} (M_V^*)$  
is away from  small neighbourhood $U$ of those intersection  points. 
Then we can choose a function $f:\R\to\R$ such that the curve 
$v=f'(u)$ is arbitrarily close to $L_{Y_1}$ and 
$L_{Y_0}$ away from the region $U$. 
This curve is illustrated in Figure \ref{fig0-tri}. 
The closeness can be measured by a small parameter $\epsilon$, 
such that as $\epsilon\to 0$, $v=f'(u)$ approaches $L_{Y_1}$ and 
$L_{Y_0}$ away from the region $U$.  We suppress the dependence of 
$v=f'(u)$  on $\epsilon $. 
 
Fix a $U(1)$-connection $A_0$ representing $(0,0)$ on $\chi(T^2)$. 
For any $U(1)$-connection $A$, define 
$T_A$ to be 
\[ 
T_z(A)= -i\int_{\{z\} \times S^1} (A-A_0), \qquad (z\in D^2). 
\] 
Choose a compactly supported 2-form $\mu$  representing the generator 
of $ H^2_{cpt} (D^2 \times S^1)$, 
such that we have 
$\int_{D^2 \times \{pt \}}\mu  = 1$ 
 for any point on $S^1$. Under the isomorphism $H^2_{cpt} (\nu(K)) 
\cong H_1(\nu(K))$, given by Poincar\'e duality, this form corresponds 
to the generator $[\mu]= PD_{\nu(K)}(l)$. 
The class of $\mu$ in $H^2( D^2 \times S^1)$ 
is trivial, and we can write $\mu= d\nu$, where $\nu$ is a 1-form 
satisfying $\int_{S^1 \times \{pt \}}\nu  = 1$, i.e. 
$\nu= PD_{T^2}(l)$. 
 
Now perturb the Chern-Simons-Dirac functional on $\nu(K) \subset Y(r)$ 
by adding the term 
\[ 
\int_{D^2} f( T_z(A) )\mu.  
\] 
Then  the perturbed Seiberg-Witten equations can be written in the following 
way: 
\ba 
\left\{\begin{array}{l} 
 F_A = *\sigma(\psi,\psi) + f'(T_A) \mu \\[2mm] 
\dirac_A (\psi ) = 0 
\end{array} 
\right. . 
\label{p-SW} 
\na 
Denote the moduli space of (\ref{p-SW}) (with generic perturbation from 
${\cal P}_0$) on $Y(r)$ by $\M_{Y, \mu}$. 
With respect to the chosen metric on $\nu(K)$, with sufficiently large 
positive scalar curvature on the support of $\mu$, the only solutions 
of the perturbed monopole equations on $\nu(K) \subset Y(r)$  
are reducibles $(A, 0)$, that satisfy 
\ba 
F_A = f'(T_A)\mu. 
\label{deformflat} 
\na 
 
With these preliminary results in place, we can prove  
the main theorem (Theorem \ref{surgery-SW3}) of this paper. 
 
\noindent{\bf Proof of Theorem \ref{surgery-SW3}.} 
This now follows from the previous discussions and the gluing 
map (cf Theorem \ref{Gluing}). From 
Theorem \ref{Gluing} and the surgery perturbation (\ref{p-SW}) on 
$\nu(K)\subset Y(r) $, we have 
$$ {\cal M}_{Y, \ \mu}^*\cong {\cal M}_{V,Y}^*\times _{\chi (T^2,Y)} 
\{ v=f'(u)\}. $$ 
 
Since we are gluing away from the lattice of $\pi^{-1}(\Theta)$, 
the limiting  points of the open ends of $\M_V^*$ and the 
neighbourhood $U$ of the intersections between the character lines, we obtain 
that solutions of the equations (\ref{p-SW}) can be identified with  
\[ 
{\cal M}_{Y, \ \mu}^*\cong {\cal M}_{V}^* 
\times _{\chi (T^2,Y)} \{ \ either \ v - u = 1, \ or \ 
u = 2k, 0\neq k\in \Z\   or \ u =\eta \}, 
\] 
when the curve $v=f'(u)$ is sufficiently close to 
the line $ \{ v - u = 1\}$ on  $\chi(T^2, Y_1)$ and the line
\[
 \{u = 2k, 0\neq k\in \Z\},\qquad  \text{ and }
\{u =\eta \}\]
 on $\chi(T^2, Y_0)$ (see Figure \ref{fig0-tri}) away 
from $U$. This shows that 
\[ 
{\cal M}_{Y, \ \mu}^*\cong  
{\cal M}_{Y_1}\cup \bigcup_{\s_k} {\cal 
M}_{Y_0}(\s_k), 
\] 
as claimed in Theorem \ref{surgery-SW3}. 
 
\section{Relative grading} 
 
In this section we show that the grading of the Floer complex 
$C_*(Y,\mu)$, defined with respect to the unique reducible point 
$\theta_Y$, induces compatible gradings on the Floer complexes 
$C_*(Y_1)$ and $C_*(Y_0,\s_k)$.  
The main tools we need in this Section are splitting formulae for 
the spectral flow, as in \cite{CLM2}, \cite{DaKi}, \cite{Nic-SF}. 
We shall first set up the necessary notation. 
 
On the space $\chi(T^2, Y_0)$ whose tangent space at any 
point  is $H^1(T^2,\R)$, we introduce the 
symplectic structure: $(a,b)\mapsto \int_{T^2} a\wedge b$, 
for $a, b \in H^1(T^2,\R)$,  
 consider the following Lagrangian 
submanifolds of $\chi(T^2, Y_0)$ 
$$\ell_{Y_1} = \pi^*(\partial_\infty({\cal M}_{\nu(K), Y_1}))= \{ (u,v)\in 
\R^2 | v-u =1 \}, $$ 
where 
$$\pi : \chi(T^2,Y_0)\to \chi(T^2,\nu(K))$$ 
is the covering map. We can identify this Lagrangian submanifold  
with a constant path of Lagrangian subspaces in 
$H^1(T^2,\R)$, given by the tangent spaces along $\ell_{Y_1} $, 
which we  denote $\tilde\ell_1(t)$. Similarly, we can consider the lines  
$$ \ell_{Y_0}(k) =\{ (2k ,v) | v\in \R \}, $$ 
for any fixed $0\neq k\in \Z$, 
and \[ 
 \ell_{Y_0}(0) =\{ (\eta, v) | v\in \R \},\] 
then we have  
$$\cup_{k\in \Z} \ell_{Y_0}(k) = 
\partial_\infty ({\cal M}_{\nu(K), Y_0}).$$ 
Each Lagrangian submanifold $\ell_{Y_0}(k)$ in $\chi(T^2,Y_0)$ determines 
a path $\tilde\ell_{Y_0}(k)$ of Lagrangian subspaces in the 
tangent space $H^1(T^2,\R)$.  
 
Moreover, there is a smooth curve 
$$\ell_\mu=\pi^*(\partial_\infty({\cal 
M}_{\nu(K), Y_1}))= \{ (u,v)\in \R^2 | v=f'(u) \} ,$$  
with $ \pi : \chi(T^2,Y_0)\to \chi(T^2,\nu(K)).$   
We can form smoothly varying 
Lagrangians of $H^1(T^2,\R)$, by taking the tangent space along the curve. 
We denote the resulting Lagrangians by $\tilde\ell_\mu$.  
 
Given any choice of two Lagrangians $\tilde\ell_\pm$ in the tangent space 
$H^1(T^2,\R)$ at the same point on $\chi(T^2, Y_0)$ we can 
define the operators that linearize  
the monopole equations 
on the manifolds with boundary $V(r)=V\cup_{T^2} 
T^2\times [0,2r]$ and $\nu(K)(r)=\nu(K)\cup_{T^2}T^2\times 
[0,2r]$. More precisely, for a sufficiently large $r\geq r_0$,  the 
gluing theorem gives  
a splitting  $(A,\psi)=(A',\psi')\#_r (a,0)$, and we can 
consider the operators (Cf. Section 4.2) on the extended $L_1^2$ spaces 
\[\begin{array}{c} 
 H_{(A',\psi'),\tilde\ell_+}: L^2_1(P_+\oplus \tilde\ell_+)\to L^2 \\[2mm] 
 H_{(A',\psi'),\tilde\ell_-}: L^2_1(P_-\oplus \tilde\ell_-)\to L^2, 
\end{array} \] 
where $P_\pm$ are APS boundary conditions \cite{APS} on the extended $L_1^2$ 
forms and spinors. 
 
Suppose we are given a path $\tilde\ell(\tau)$ of Lagrangians  
in $H^1(T^2,\R)$, which can be written in the form 
$\tilde\ell(\tau)=T_{a(\tau)}\ell$,  
for some Lagrangian submanifold $\ell$ of $\chi(T^2, Y_0)$ with a 
regular parameterization $a(\tau)$. Assume that, for $0\leq \tau \leq 1$ 
the arc $a(\tau)$ of the Lagrangian submanifold $\ell$ avoids the lattice 
of $\{ \pi^{-1} (\Theta) \}$ and the limiting points of 
$\d_\infty (M^*_V)$ on the circle 
$\chi(V)$.  
Moreover, we  
assume that we have $a$ and $b$ in $\ell\cap \ell_{Y_1}$ and that 
$\ell$ and $\ell_{Y_1}$ intersect transversely. Assume the arc of $\ell_{Y_1}$ 
between these endpoints is parameterized over the same interval $0\leq 
\tau \leq 1$. Moreover, for small enough $\epsilon$ in the surgery 
perturbation $\mu$, there are  
distinct points $a^\epsilon$ and $b^\epsilon$ 
in $\ell\cap \ell_\mu$. We can assume that, for $\epsilon$ 
sufficiently small, also $\ell$ and $\ell_\mu$ intersect transversely, 
and there are parameterizations of the arcs of Lagrangians $\ell$ and 
$\ell_\mu$ with 
endpoints $a^\epsilon$ and $b^\epsilon$. 
 
We have the following result, which is the key lemma in the 
comparison of the Maslov indices. 
 
\begin{Lem} 
With the hypothesis as above, we have 
$$ Maslov(\tilde\ell(\tau),\tilde \ell_{Y_1})= 
Maslov(\tilde\ell(\tau),\tilde\ell_\mu(\tau)), $$ 
where the first Maslov index is computed with respect to the 
parameterizations with endpoints $a$ and $b$, and the second 
with respect to the parameterizations with endpoints $a^\epsilon$ 
and $b^\epsilon$, as specified above.  
\label{equal-maslov-1-mu} 
\end{Lem} 
 
\proof By applying the properties of the Maslov index 
(cf. \cite{CLM-Maslov}, Section 1), we can see that the claim follows, 
upon showing that we have  
$$  Maslov(\tilde\ell_\mu(\tau), \tilde\ell_{Y_1})=0 $$ 
which is obvious by the choice of $\ell_\mu(\tau)$ and $\ell_{Y_1}$. 
\endproof 
 
As a consequence of this result, we obtain the following proposition 
relating the relative gradings on $\M_{Y,\mu}^*$ and $\M^*_{Y_1}$ 
respectively. Given a path $\{ (A'(\tau),\psi'(\tau)) | \tau\in [0,1]
\}$, and a corresponding path $\{ (a(\tau),0) | \tau\in [0,1] \}$, we
can consider the
corresponding paths of operators $H_{(A'(\tau),\psi'(\tau))}$,
$H_{(a(\tau),0)}$, and $H_{(A'(\tau),\psi'(\tau))\#_r(a(\tau),0)}$.

\begin{Pro} 
Suppose we are given two irreducible critical points $a=[A_a,\psi_a]$ and 
$b=[A_b,\psi_b]$ in $\M^*_{Y_1}$, and the corresponding elements  
$a^\epsilon=[A_a^\epsilon,\psi_a^\epsilon]$ and 
$b^\epsilon=[A_b^\epsilon,\psi_b^\epsilon]$ in $\M_{Y,\mu}^*$. Then we have 
$$ \deg_{Y,\mu}(a^\epsilon)-\deg_{Y,\mu}(b^\epsilon)= 
\deg_{Y_1}(a)-\deg_{Y_1}(b). $$ 
\label{same-deg-1-mu} 
\end{Pro} 
\proof Under the pre-gluing map, we can assume that 
$(A_a^\epsilon,\psi_a^\epsilon)$ and $(A_b^\epsilon,\psi_b^\epsilon)$ are 
connected by a path $(A'(\tau),\psi'(\tau))\#_r(a(\tau),0)$  
($\tau \in [0, 1]$), where we have 
$$ (A_a^\epsilon,\psi_a^\epsilon) 
=(A'(0),\psi'(0))\#_r(a(0),0), $$ 
$$ (A_b^\epsilon,\psi_b^\epsilon) = 
(A'(1),\psi'(1))\#_r(a(1),0). $$ 
Then by definition,  
$$ \deg_{Y,\mu}(A_a^\epsilon,\psi_a^\epsilon)- 
\deg_{Y,\mu}(A_b^\epsilon,\psi_b^\epsilon) =  
 \frac{1}{r^2} SF_{Y(r)}( H_{(A'(\tau),\psi'(\tau))\#_r(a(\tau),0)}), $$ 
We can compute this spectral flow with the splitting formula on $Y(r)$ 
from (Theorem C of \cite{CLM2}). 
We obtain 
$$ \epsilon SF (H_{(A'(\tau),\psi'(\tau)),\tilde\ell(\tau)}) +  Maslov 
(\tilde\ell(\tau),\tilde\ell_\mu) + \epsilon SF 
(H_{(a(\tau),0),\tilde\ell_\mu}). $$  
With the analogous splitting 
formula on $Y_1(r)$, by applying 
the Capell-Lee-Miller decomposition of the spectral 
flow (Theorem C of \cite{CLM2}), we obtain 
\[ 
\begin{array}{lll}&& \deg_{Y_1}(A_a,\psi_a) - \deg_{Y_1}(A_b,\psi_b)\\[2mm] 
&=& \epsilon SF (H_{(A'(\tau),\psi'(\tau)),\tilde\ell(\tau)}) +  Maslov 
(\tilde\ell(\tau),\tilde\ell_1) + \epsilon SF 
(H_{(a(\tau),0),\tilde\ell_1}).  
\end{array} 
\]  
In both cases, we can assume that we consider the 
same boundary value problem (the same choice of Lagrangians) for the 
operator on the knot complement $V$. 
We choose $\tilde\ell_\mu$ or $\tilde\ell_1$ for the operator on the 
tubular neighbourhood of the knot $\nu(K)$. The previous Lemma 
shows that the quantities $\epsilon SF 
(H_{(a(\tau),0),\tilde\ell_1})$ and $\epsilon SF 
(H_{(a(\tau),0),\tilde\ell_{\mu}})$ coincide, and that the two Maslov 
indices are also the same. 
\endproof 
 
Similarly, we can now compare the relative grading of two solutions in 
$\M_{Y_0}(\s_k)$ with the relative grading of the corresponding 
solutions in $\M_{Y,\mu}$. 
 
Again, suppose we are given a path $\tilde\ell(\tau)$  
of Lagrangians in the tangent 
space $H^1(T^2,\R)$, of the form $\tilde\ell(\tau)=T_{a(\tau)}\ell$, 
for some Lagrangian submanifold $\ell$ of $\chi(T^2,V)$ with a 
regular parameterization $a(\tau)$. Assume that, for $0\leq \tau \leq 1$ 
the arc $a(\tau)$ of the Lagrangian submanifold $\ell$ avoids the lattice 
of $\{ \pi^{-1}(\Theta) \}$ and the limiting points  
$\d_\infty(\M_V^*)$ on the circle $\chi(V)$. Moreover, we  
assume that we have $a$ and $b$ in $\ell\cap \ell_{Y_0}(k)$ and that 
$\ell$ and $\ell_{Y_0}(k)$ intersect transversely.  
Assume the arc of $\ell_{Y_0}(k)$ 
between these endpoints is parameterized over the same interval $0\leq 
\tau \leq 1$. Moreover, for small enough $\epsilon$ in the surgery 
perturbation $\mu$, there are points $a^\epsilon$ and $b^\epsilon$ 
in $\ell\cap \ell_\mu$. We can assume that, for $\epsilon$ 
sufficiently small, also $\ell$ and $\ell_\mu$ intersect transversely, 
and there are parameterizations of the arcs of Lagrangians $\ell$ and 
$\ell_\mu$ with endpoints $a^\epsilon$ and $b^\epsilon$. 
 
With these hypothesis we have the following lemma, whose 
proof is analogous to the proof of Lemma \ref{equal-maslov-1-mu}. 
 
\begin{Lem} 
With the hypothesis as above,  we have 
$$ Maslov(\tilde\ell(\tau),\tilde\ell_{Y_0}(k))= 
Maslov(\tilde\ell(\tau),\tilde\ell_\mu(\tau)), $$ 
where the first Maslov index is computed with respect to the 
parameterizations with endpoints $a_0$ and $b_0$, and the second 
with respect to the parameterizations with endpoints $a_0^\epsilon$ 
and $b_0^\epsilon$, as specified above.  
\label{equal-maslov-0-mu} 
\end{Lem} 
 
We have the following proposition relating the relative gradings 
on $\M_{Y,\mu}^*$ and ${\cal M}_{Y_0,\s_k}$ (for $k\in \Z$) respectively. 
 
\begin{Pro} 
Suppose we are given $a=[A_a,\psi_a]$ and $b=[A_b,\psi_b]$ representing two 
elements in ${\cal M}_{Y_0,\s_k}$, and let 
$a^\epsilon=[A_a^\epsilon,\psi_a^\epsilon]$ and 
$b^\epsilon=[A_b^\epsilon,\psi_b^\epsilon]$ be the corresponding elements in 
$\M_{Y,\mu}^*$. We have  
$$ \deg_{Y_0,\s_k}(A_a,\psi_a)- 
\deg_{Y_0,\s_k}(A_b,\psi_b) = \deg_{Y,\mu}(A_a^\epsilon,\psi_a^\epsilon)- 
\deg_{Y,\mu}(A_b^\epsilon,\psi_b^\epsilon) \ \ \hbox{ mod } 2k. $$ 
\label{same-deg-0-mu} 
\end{Pro} 
\proof 
With the notation as in the Lemma \ref{equal-maslov-0-mu}, we have  
\[
(A_a,\psi_a)=(A'(0),\psi'(0))\# (a(0),0);
\]
 and 
\[
(A_b,\psi_b)= (A'(1),\psi'(1))\# (a(1),0).
\]
 We can calculate the 
relative grading using the splitting formula on $Y_0(r)$ 
$$ \begin{array}{lll}&&\deg_{Y_0,\s_k}(A_a,\psi_a)- 
\deg_{Y_0,\s_k}(A_b,\psi_b)  \\[2mm] 
&=&(\epsilon SF)(H_{(A'(\tau),\psi'(\tau)),\tilde\ell(\tau)}) + 
Maslov(\tilde\ell(\tau),\tilde\ell_0(\tau)) + 
\epsilon SF(H_{(a(\tau),0),\tilde\ell_0(\tau)})\end{array}. $$  
We can then compare directly these terms with the corresponding terms 
in the splitting formula for the spectral flow of the operators on 
$Y(r)$, as in the case of Corollary \ref{same-deg-1-mu}. The result of 
Lemma \ref{equal-maslov-0-mu} guarantees that we obtain the same 
result.  
\endproof 
 
Notice that the results of Lemma \ref{equal-maslov-0-mu} and Corollary 
\ref{same-deg-0-mu} imply that the grading $\deg_{Y,\mu}$ defines a 
choice of an integer lift of the $\Z_{2k}$-valued relative grading of 
$C_*(Y_0,\s_k)$ given by 
$$\deg_{Y_0,\s_k}(A_a,\psi_a)- 
\deg_{Y_0,\s_k}(A_b,\psi_b) = \deg_{Y,\mu}(A_a,\psi_a)- 
\deg_{Y,\mu}(A_b,\psi_b), $$ 
under the identification ${\cal M}_{Y_0,\s_k}\hookrightarrow {\cal 
M}_{Y,\mu}$. We will discuss the properties of the 
integer lift $C_{(*)}(Y_0,\s_k)$ of $C_*(Y_0,\s_k)$  elsewhere.

\vskip .2in

\noindent {\bf A.L. Carey}, 

\noindent School of Mathematical Sciences,
Australian National University, Canberra ACT, Australia\par
\noindent acarey@wintermute.anu.edu.au\par
 
\vskip .2in 

\noindent {\bf M. Marcolli}, 

\noindent Max-Planck-Institut f\"ur Mathematik, D-53111 Bonn, Germany 
\par  
\noindent marcolli\@@mpim-bonn.mpg.de 

\vskip .2in 
 
\noindent {\bf B.L. Wang}, 

\noindent Department of Pure 
Mathematics, University of Adelaide, Adelaide SA 5005 \par 
\noindent bwang\@@maths.adelaide.edu.au \par

\end{document}